\documentclass[11pt]{article}

\usepackage[T1]{fontenc}
\usepackage[a4paper,margin=27mm]{geometry}
\usepackage{amsmath,amssymb,amsthm,mathtools,mathrsfs}
\usepackage{microtype}
\usepackage[colorlinks=true,linkcolor=blue,citecolor=blue,urlcolor=blue]{hyperref}
\hypersetup{
  pdftitle={Characteristic Chow Polynomials of Layered-Permutation Intervals},
  pdfauthor={Jiahao Zhang}
}

\newtheorem{theorem}{Theorem}[section]
\newtheorem{proposition}[theorem]{Proposition}
\newtheorem{lemma}[theorem]{Lemma}
\newtheorem{corollary}[theorem]{Corollary}
\newtheorem{conjecture}[theorem]{Conjecture}
\newtheorem*{remark}{Remark}

\allowdisplaybreaks

\title{Characteristic Chow Polynomials of Layered-Permutation Intervals}
\author{Jiahao Zhang\\[4pt]
\small College of Mathematical Sciences\\
\small Tianjin Normal University\\
\small Tianjin 300387, P. R. China\\[2pt]
\small \href{mailto:jiahaozhang1215@163.com}{\texttt{jiahaozhang1215@163.com}}}
\date{}

\begin{document}

\maketitle

\begin{abstract}
We study characteristic Chow polynomials of principal intervals in
generalized subword order, equivalently of pattern-containment intervals
below layered permutations.  A terminal-chain recurrence and canonical
rightmost embeddings yield a recursive rational generating function for
every fixed number of layers, governed by a universal denominator.  In the
two-layer case we obtain an explicit kernel, prove canonical interlacing
when one layer has size at most four, and prove eventual interlacing at
every fixed width by a Laguerre limit.  We also treat two families with
unboundedly many layers.  Singleton-prefix intervals have a closed
Delannoy--Jacobi gamma polynomial and form real-rooted interlacing
sequences.  Singleton-sandwich intervals admit rational three-variable series and a
marked-slot formula for all gamma coefficients. Their central-layer-two
and central-layer-three slices are real-rooted, with interlacing under
growth from two to three. For each fixed pair of singleton flanks,
a Laguerre differential transform proves eventual interlacing as the
central layer grows.  These results exhibit composition
intervals as a structured testing ground for positivity and root-location
questions for Chow polynomials of posets.
\end{abstract}

\medskip
\noindent\textbf{Keywords.}
Characteristic Chow polynomial; composition poset; generalized subword
order; layered permutation; gamma-positivity; real-rootedness; interlacing.

\smallskip
\noindent\textbf{2020 Mathematics Subject Classification.}
Primary 06A07; Secondary 05A05, 05A15, 26C10.

\section{Introduction}
\label{sec:introduction}

Characteristic Chow polynomials form a recently developed family of
invariants for finite graded posets.  They arise from an incidence-algebra
construction parallel to Kazhdan--Lusztig--Stanley theory and, for the
lattice of flats of a matroid, recover the Hilbert--Poincar\'e series of its
Chow ring \cite{FerroniMatherneVecchi2026}.  The Hodge theory of matroids
therefore supplies their original geometric model
\cite{AdiprasitoHuhKatz2018}.  In the broader poset setting, characteristic
Chow polynomials remain palindromic, nonnegative, and unimodal; they are
gamma-positive for Cohen--Macaulay posets.  Stump obtained explicit
gamma-positive expansions for posets admitting an R-labeling
\cite{Stump2025}.  Real-rootedness is subtler: even Cohen--Macaulayness
does not suffice, as the distributive-lattice counterexample in
\cite[Remark~4.26]{FerroniMatherneVecchi2026} shows.  It has recently been established for
uniform matroids and maximal ranked posets
\cite{BrandenVecchi2026Uniform}, for suitable simplicial posets
\cite{HosterStump2025}, and for broad families controlled by totally
nonnegative matrices \cite{BrandenVecchi2025TN}, as well as for
uniformly monotonically edge-lexicographically shellable posets
\cite{CoronFerroniLi2025}.

This paper studies characteristic Chow polynomials for a different and
very concrete family of posets.  Let \(\mathcal C\) be the set of
compositions ordered by generalized subword order, and let
\(P_\beta=[\varnothing,\beta]\) be the principal interval generated by a
composition \(\beta\).  Generalized subword order on a chain was studied
through its M\"obius function and topology in
\cite{SaganVatter2006,McNamaraSagan2012}.  It also has a direct
permutation-pattern interpretation.  Writing a composition
\(\beta=(b_1,\ldots,b_k)\) as the direct sum of decreasing permutations of
lengths \(b_1,\ldots,b_k\) identifies \(P_\beta\) with the pattern-containment
interval below the corresponding layered permutation
\cite{McNamaraSteingrimsson2015}.  Thus the number of parts of \(\beta\) is
the number of layers of the host permutation.

Our starting point is a terminal-chain recurrence together with an exact
description of what changes when the final layer grows.  For two layers this
leads to the rational kernel in Theorem~\ref{thm:two-layer-kernel}.  The same
kernel supports two complementary kinds of information.  At small width it
permits exact Sturm arguments: Theorem~\ref{thm:min-four-interlacing} proves
canonical interlacing whenever one layer has size at most four.  At every
fixed width it also admits an exact transfer formula and a Laguerre limit;
Theorem~\ref{thm:eventual-fixed-width} consequently proves canonical
interlacing for all sufficiently long second layers.  These results isolate
the unsolved uniform boundary as the explicit gamma-difference problem in
Conjecture~\ref{conj:gamma-difference-sturm}.  The use of interlacing,
orthogonal-polynomial limits, and gamma expansions follows a productive
line of work on Eulerian- and Narayana-type polynomials; see, for example,
\cite{ChenYangZhang2018,HaglundZhang2019,LinMaZhang2021,YangZhang2025,Zhang2019Veronese}.

The two-layer calculation is not the endpoint of the paper.  Its
rightmost-embedding decomposition persists in every fixed dimension.
Theorem~\ref{thm:fixed-layer-rationality} gives, for each fixed number
\(k\) of layers, a recursive rational generating function for all principal
intervals with \(k\) parts.  The common transfer denominator is
\[
 D_k(\mathbf x;t)=
 \sum_{m=0}^{k}(-1)^m[m+1]_t e_m(\mathbf x),
\]
which specializes to the two-layer kernel and produces its three-layer
analogue without a separate case analysis.  This denominator also occurs
in the Toeplitz-matrix Chow series of \cite{BrandenVecchi2025TN}; the
composition-specific content is the canonical embedding operator and the
recursive numerators, not a new denominator.

We then pass beyond fixed dimension. The nonempty part of the interval below
\((1^r,b)\), whose number of layers is unbounded, is a product of two
chains. Theorem~\ref{thm:singleton-prefix-rectangles} identifies its gamma
polynomial with a Delannoy--Jacobi polynomial and proves real-rootedness and
consecutive interlacing for all \(r\) and \(b\). Allowing singleton layers
on both sides of the unique nonsingleton layer gives the three-parameter
family \((1^r,b,1^s)\). Theorem~\ref{thm:singleton-sandwich-series} derives
complete rational Chow and gamma series for this family, and
Proposition~\ref{prop:sandwich-slot-model} interprets every gamma coefficient
through selected flank positions and marked central slots.
Theorems~\ref{thm:central-two-real-rootedness}
and~\ref{thm:central-three-interlacing} prove real-rootedness for all flank
lengths when the central layer has size two or three, including interlacing
under growth from two to three. The proof combines a Jacobi integral
transform with an explicit sign comparison at its zeros.
Finally, Theorem~\ref{thm:eventual-sandwich} proves eventual real-rootedness
and consecutive interlacing as the central layer grows, for every fixed
pair of flanks, through a Laguerre differential transform.

The main point is therefore structural rather than computational.
Canonical embeddings resolve the multiple-containment issue at each fixed
layer number, while the singleton families turn interval structure into
explicit gamma coefficients and parameter-uniform root comparisons.
The Jacobi and Laguerre arguments establish the stated ranges without
assuming real-rootedness for all composition intervals.  Together these
results connect generalized subword order and layered permutation patterns
with Chow-polynomial root-location questions beyond the standard
rank-uniform setting.

The remainder of the paper follows this progression.  Section
\ref{sec:composition-two-layer} develops the terminal recurrence and the
two-layer theory.  The fixed-layer construction appears next, followed by
the two unbounded-layer families.  The final section records the resulting
open problems and explains the role of the accompanying exact checks.


\section{Composition intervals and the two-layer kernel}
\label{sec:composition-two-layer}

Let \(\mathcal C\) be the composition poset under generalized subword order.
Thus \(\alpha=(a_1,\ldots,a_r)\leq\beta=(b_1,\ldots,b_k)\) if and only if
there are indices \(i_1<\cdots<i_r\) such that \(a_j\leq b_{i_j}\) for every
\(j\); see \cite{SaganVatter2006,McNamaraSagan2012}.  If
\(\delta_m=m(m-1)\cdots1\) and
\[
 \pi_\beta=\delta_{b_1}\oplus\cdots\oplus\delta_{b_k},
\]
then \(\pi_\beta\) is layered, and \(\beta\mapsto\pi_\beta\) identifies
generalized subword order with pattern containment among layered
permutations \cite{McNamaraSteingrimsson2015}.  In particular, the principal
interval below \(\pi_\beta\) is isomorphic to
\(P_\beta=[\varnothing,\beta]\).  We rank this interval by total size and
denote its characteristic Chow polynomial by \(H_\beta(t)\), with
\(H_\varnothing(t)=1\).  As usual,
\([d]_t=1+t+\cdots+t^{d-1}\) for \(d\geq1\), while \([d]_t=0\) for
\(d\leq0\).

To fix the convention, for a graded bounded poset $P$ let
\[
 \chi_{xy}(t)=\sum_{x\leq u\leq y}\mu(x,u)t^{\rho(y)-\rho(u)},
 \qquad
 \overline\chi_{xy}(t)=
 \begin{cases}
  \chi_{xy}(t)/(t-1),&x<y,\\
  -1,&x=y.
 \end{cases}
\]
Here $\mu$ is the M\"obius function and $\rho$ is the rank function.
The characteristic Chow function is $H=-\overline\chi^{-1}$, where
inversion uses incidence-algebra convolution, and $H_\beta$ is its
$(\varnothing,\beta)$ entry; see \cite[Definition~3.3]{FerroniMatherneVecchi2026}.
We use the equivalent chain formula below.

For a nonempty composition \(\alpha\), define its terminal-chain polynomial
by
\[
 R_\alpha(t)=
 \sum_{\varnothing=\alpha_0<\cdots<\alpha_m=\alpha}
 t^{m-1}\prod_{i=1}^m
 [|\alpha_i|-|\alpha_{i-1}|-1]_t.
\]
Grouping the Chow chain formula of
\cite[Theorem~4.1]{FerroniMatherneVecchi2026} by the terminal element gives
\begin{equation}
 H_\beta(t)=1+t\sum_{\varnothing<\alpha\leq\beta}R_\alpha(t),
 \label{eq:terminal-sum}
\end{equation}
and grouping a terminal chain by its penultimate element gives
\begin{equation}
 R_\alpha(t)=[|\alpha|-1]_t+
 t\sum_{\varnothing<\nu<\alpha}
 R_\nu(t)[|\alpha|-|\nu|-1]_t.
 \label{eq:terminal-recursion}
\end{equation}

The chain formula above is the specialization to \(P_\beta\) of the
characteristic Chow function of a graded bounded poset
\cite{FerroniMatherneVecchi2026}. For nonempty $\beta$, \(H_\beta(t)\) is
palindromic of degree \(|\beta|-1\), and hence has a unique gamma expansion
\[
 H_\beta(t)=\sum_{i=0}^{\lfloor(|\beta|-1)/2\rfloor}\gamma_i(H_\beta)
 t^i(1+t)^{|\beta|-1-2i}.
\]
We call \(\Gamma_\beta(s)=\sum_i\gamma_i(H_\beta)s^i\) its gamma
polynomial, taking coefficients outside the displayed range to be zero.
A palindromic polynomial is \emph{gamma-positive} if these
coefficients are nonnegative.

For real-rooted polynomials of degrees differing by at most one, we write
\(f\preceq g\) when their zeros weakly interlace, and \(f\prec g\) when the
interlacing is strict. Here these symbols record alternation, without
choosing an orientation when the degrees are equal; constant--linear pairs
are included. We repeatedly use the elementary gamma--Chow root
correspondence.  Under
\[
 \mathcal T_n p(t)=(1+t)^n
 p\!\left(\frac{t}{(1+t)^2}\right),
\]
each negative zero \(s\) of \(p\) gives the reciprocal pair
\[
 \frac{1-2s+\sqrt{1-4s}}{2s},
 \qquad
 \frac{1-2s-\sqrt{1-4s}}{2s},
\]
both on the negative real axis, while unused degree contributes roots at
\(-1\).  The precise root orders needed for transferring interlacing are recorded
in Lemma~\ref{lem:gamma-chow-interlacing} below; mere alternation without
its direction would not suffice.  This substitution is a standard tool in
real-rooted and gamma-positive enumeration; compare
\cite{ChenYangZhang2018,HaglundZhang2019,Zhang2019Veronese}.

We also use the one-variable convex-cone property of proper position
\cite{BorceaBranden2010}.  In the instances below the compatible
orientation is established explicitly by signs at consecutive zeros.  If
\(f\) and \(g\) have that alternating order and positive leading
coefficients, then evaluation at the zeros of either summand shows that
every positive linear combination of \(f\) and \(g\) remains in the
same interlacing cone.

\begin{lemma}[Gamma-to-Chow interlacing]
\label{lem:gamma-chow-interlacing}
Let $p$ and $q$ have positive constant terms and only simple negative zeros,
with $d=\deg p$, $e=\deg q$, $2d\leq n$, and $2e\leq n+1$.
Write their zeros increasingly as $\alpha_i$ and $\beta_i$, respectively.
Suppose either
\[
 e=d,\qquad
 \alpha_1<\beta_1<\alpha_2<\beta_2<\cdots<\alpha_d<\beta_d,
\]
or
\[
 e=d+1,\qquad
 \beta_1<\alpha_1<\beta_2<\cdots<\alpha_d<\beta_{d+1}.
\]
Empty lists are allowed. Then $\mathcal T_np\preceq\mathcal T_{n+1}q$.
Their only possible common zero is $-1$.
\end{lemma}

\begin{proof}
If $d=e=0$, both transforms are positive multiples of powers of
$1+t$, and the assertions are immediate. Otherwise, proceed as follows.
If $a$ is the leading coefficient of $p$, then
\[
 \mathcal T_np(t)=a(1+t)^{n-2d}
       \prod_{i=1}^d\bigl(t-\alpha_i(1+t)^2\bigr).
\]
Thus its degree is exactly $n$, its multiplicity at $-1$ is exactly
$m=n-2d$, and every other zero belongs to a reciprocal pair, one on
each side of $-1$. The same statements hold for $q$.
The map $t\mapsto t/(1+t)^2$ is strictly decreasing on
$(-\infty,-1)$ and strictly increasing on $(-1,0)$, onto the negative
half-line in both cases. Consequently, the order of the gamma zeros is
reversed on the left branch and preserved on the right branch.

When $e=d$, the left-branch alternation ends with a zero of
$\mathcal T_np$, and the right-branch alternation starts with one.
The multiplicities at $-1$ are $m$ and $m+1$, so the intervening lists
alternate starting and ending with $\mathcal T_{n+1}q$.
When $e=d+1$, both branch lists meet $-1$ at a zero of
$\mathcal T_{n+1}q$. Here $m\geq1$, and the multiplicities are $m$
and $m-1$, so the intervening lists alternate starting and ending with
$\mathcal T_np$. This proves the full weak alternation, including
empty branch lists. Strictness of the gamma orders and injectivity of
each branch exclude any common zero other than $-1$.
\end{proof}

\begin{lemma}[Interlacing of a Chow increment]
\label{lem:gamma-increment-transfer}
Let $p(0)>0$, and let $g\neq0$ have positive leading coefficient and
$g(0)=0$. Suppose that $p$ has only negative zeros, $g$ has only
nonpositive zeros, and $p\preceq g$. Assume $2\deg p\leq n$ and
$2\deg g\leq n+1$. Put
\[
 A=\mathcal T_np,\qquad B=\mathcal T_{n+1}(p+g),\qquad C=B-A.
\]
Then $A,B,C$ are real-rooted and $A\preceq C\preceq B$.
\end{lemma}

\begin{proof}
First suppose the interlacing is strict. Since $0$ is the largest zero
of $g$, its degree is either $\deg p$ or $\deg p+1$, and partial fractions give
\[
 \frac{g(s)}{p(s)}=a s+\sum_i c_i\frac{s}{s-\alpha_i},
 \qquad a\geq0,\quad c_i>0,\quad \alpha_i<0.
\]
Indeed, the residues at the zeros $\alpha_i$ of $p$ are negative by
the alternating signs, and the value at $0$ determines the constant term.
Consequently,
\[
 \frac{C(t)}{A(t)}
 =t+a\frac{t}{1+t}
   +\sum_i c_i\frac{t(1+t)}{t-\alpha_i(1+t)^2}.
\]
After cancelling common factors, at a pole $\tau\neq-1$ the corresponding residue is
$c_i\tau(1+\tau)^2/(1-\tau)<0$; at $-1$ the residue is $-a$.
The pole at $-1$ occurs only when $a>0$; then $n-2\deg p\geq1$,
so the displayed expression still represents a polynomial quotient.
The polynomial part has slope $1$. Writing the quotient as
$t+\kappa+\sum_\tau r_\tau/(t-\tau)$, with $r_\tau<0$, shows that
its derivative is $1-\sum_\tau r_\tau/(t-\tau)^2>0$.
Thus $C/A$ increases from $-\infty$ to $+\infty$ on every interval
between poles, including the two exterior intervals.
The level sets $C/A=0$ and $C/A=-1$ give alternating real zeros of
$C$ and $B=C+A$, interlaced by the poles. If $M$ is the number of
poles, their reduced numerators both have degree $M+1$, so the level
sets exhaust all zeros. The common factor to restore is
$(1+t)^{n-2\deg p}$ when $a=0$, and
$(1+t)^{n-2\deg p-1}$ when $a>0$; inserting its identical root
lists preserves weak alternation. For weak
interlacing, approximate the ordered zeros by strictly alternating
negative zeros, keeping the largest zero of $g$ at $0$ and the leading
coefficients positive. The transforms have fixed degrees and nonzero
limiting leading coefficients, so real-rootedness and weak alternation
persist in the limit.
\end{proof}


\subsection{The exact growth filter}

For a nonempty composition \(\beta=(b_1,\ldots,b_k)\), set
\[
 \beta^-=
 \begin{cases}
 (b_1,\ldots,b_k-1),&b_k>1,\\
 (b_1,\ldots,b_{k-1}),&b_k=1.
 \end{cases}
\]

\begin{lemma}[Final-layer filter]
\label{lem:final-layer-filter}
Let \(\beta=(u,b)\).  Then
\[
 P_\beta\setminus P_{\beta^-}
 =\{(v,b):v\leq u\text{ and }(v,b)\not\leq u\}.
\]
If \(j\) is the last position of \(u\) whose part is at least \(b\), and
\(q\) is the prefix of \(u\) strictly before position \(j\), deletion of the
final part is an order isomorphism
\[
 P_\beta\setminus P_{\beta^-}\cong P_u\setminus P_q.
\]
If there is no such position, deletion of the final part gives an order
isomorphism with all of \(P_u\).
\end{lemma}

\begin{proof}
The old interval is a lower ideal.  Any new composition must therefore use
the final coordinate \(b\) in every embedding into \((u,b)\).  When \(b>1\),
the part placed there must equal \(b\), since a smaller part would still fit
after replacing \(b\) by \(b-1\); when \(b=1\), positivity gives the same
conclusion.  Hence a new element has the form \((v,b)\), where \(v\leq u\),
and it is new exactly when \((v,b)\not\leq u\).

Suppose that position \(j\) is the last position of \(u\) capable of holding
\(b\).  Placing the final \(b\) as far right as possible shows that
\((v,b)\leq u\) if and only if \(v\leq q\).  If there is no capable position,
\((v,b)\not\leq u\) for every \(v\leq u\).  Finally,
\((v,b)\leq(w,b)\) if and only if \(v\leq w\), which proves the asserted
order isomorphisms.
\end{proof}

In particular, the new part is an upper filter, and a chain crosses from the
old interval to the new part at most once.  Equations
\eqref{eq:terminal-sum} and Lemma~\ref{lem:final-layer-filter} give the exact
growth identity
\begin{equation}
 H_{(u,b)}(t)-H_{(u,b)^-}(t)
 =t\sum_{\substack{v\leq u\\(v,b)\not\leq u}}R_{(v,b)}(t).
 \label{eq:final-layer-growth}
\end{equation}

The first gamma coefficient already has a simple interpretation.

\begin{proposition}
\label{prop:first-gamma}
If \(|\beta|\geq2\), then
\[
 \gamma_1(H_\beta)=|P_\beta|-|\beta|-1.
\]
Consequently,
\[
 \gamma_1(H_\beta)-\gamma_1(H_{\beta^-})
 =|P_\beta\setminus P_{\beta^-}|-1.
\]
\end{proposition}

\begin{proof}
Only chains with one nontrivial step contribute to \([t]H_\beta(t)\), and
their terminal element may be any element of rank at least two.  The unique
rank-one element is \((1)\), so \([t]H_\beta(t)=|P_\beta|-2\).  Comparing this
with the coefficient of \(t\) in
\[
 H_\beta(t)=\sum_i\gamma_i(H_\beta)
 t^i(1+t)^{|\beta|-1-2i}
\]
proves the first identity.  The first formula also gives zero for the rank-one interval, taking
$\gamma_1=0$ for degrees below two. Thus the second follows by subtraction,
including $|\beta|=2$.
\end{proof}

\subsection{A rational kernel for all two-layer intervals}

For a one-part composition, recurrence \eqref{eq:terminal-recursion} gives
\begin{equation}
 R_{(s)}(t)=(1+t)^{s-2}\quad(s\geq2),
 \qquad R_{(1)}(t)=0.
 \label{eq:one-part-terminal}
\end{equation}
Let
\[
 \mathcal R(x,y;t)=\sum_{a,b\geq1}R_{(a,b)}(t)x^ay^b.
\]

\begin{theorem}[Two-layer rational kernel]
\label{thm:two-layer-kernel}
Put \(c=1+t\).  Then
\begin{align}
 \mathcal R(x,y;t)
 &=\frac{A(x,y;t)}{1-tK(x,y;t)},
 \label{eq:R-kernel}\\
 A(x,y;t)
 &=\frac{xy(1-c^2xy)}
 {(1-cx)(1-cy)(1-cxy)(1-tcxy)},
 \label{eq:A-kernel}\\
 K(x,y;t)
 &=\frac{x^2+xy+y^2-cxy(x+y)+tx^2y^2}
 {(1-x)(1-y)(1-tx)(1-ty)}.
 \label{eq:K-kernel}
\end{align}
Moreover,
\begin{equation}
 1-tK(x,y;t)=
 \frac{1-c(x+y)+(1+t+t^2)xy}
 {(1-x)(1-y)(1-tx)(1-ty)},
 \label{eq:kernel-cancellation}
\end{equation}
and
\begin{equation}
 H_{(a,b)}(t)=c^{\max(a,b)-1}
 +t\sum_{i=1}^a\sum_{j=1}^bR_{(i,j)}(t).
 \label{eq:H-from-R}
\end{equation}
\end{theorem}

\begin{proof}
The elements below \((a,b)\) are the empty composition, the one-part
compositions \((s)\) with \(s\leq\max(a,b)\), and the two-part compositions
\((i,j)\) with \(i\leq a\), \(j\leq b\).  By
\eqref{eq:one-part-terminal}, the contribution of the empty composition and
the one-part compositions to \eqref{eq:terminal-recursion} is
\[
 A_{a,b}(t)=c^{\max(a,b)-1}[\min(a,b)]_t.
\]
The remaining terms form a bivariate convolution with kernel
\([r+s-1]_t\).  Hence, with
\[
 A(x,y;t)=\sum_{a,b\geq1}A_{a,b}(t)x^ay^b,
 \qquad
 K(x,y;t)=\sum_{r,s\geq0}[r+s-1]_t x^ry^s,
\]
recurrence \eqref{eq:terminal-recursion} becomes
\(\mathcal R=A+t\mathcal RK\), which proves \eqref{eq:R-kernel}.

Splitting the sum defining \(A\) into \(a\geq b\) and \(b>a\), and summing
the two geometric series, gives \eqref{eq:A-kernel}.  Summing the
\(t\)-integers in \(K\) gives \eqref{eq:K-kernel}.  On placing the latter
over a common denominator, all terms of total degree at least two cancel
except the mixed term, giving \eqref{eq:kernel-cancellation}.  Finally,
\eqref{eq:terminal-sum} and
\[
 1+t\sum_{s=2}^{\max(a,b)}c^{s-2}=c^{\max(a,b)-1}
\]
give \eqref{eq:H-from-R}.
\end{proof}

The kernel becomes more transparent after passing to gamma polynomials.  Set
\[
 H_{(a,b)}(t)=c^{a+b-1}
 \Gamma_{a,b}\!\left(\frac{t}{c^2}\right),
 \qquad
 \mathcal G(X,Y;s)=\sum_{a,b\geq1}\Gamma_{a,b}(s)X^aY^b.
\]

\begin{corollary}[Two-layer gamma kernel]
\label{cor:two-layer-gamma-kernel}
We have
\[
 \mathcal G(X,Y;s)=
 \frac{XY(1-XY)N(X,Y;s)}
 {(1-X)(1-Y)(1-XY+sX^2Y^2)
  \bigl(1-X-Y+(1-s)XY\bigr)},
\]
where
\[
 N(X,Y;s)=1-(1-s)(X+Y)+(1-3s)XY
 +sXY(X+Y)-s(1-s)X^2Y^2.
\]
Consequently, $\deg\Gamma_{a,b}\leq\min(a,b)$.
In particular, the non-diagonal factor is the Delannoy--Jacobi kernel
\[
 \frac{1}{1-X-Y+(1-s)XY}
 =\sum_{p,q\geq0}\sum_k
 \binom pk\binom qk s^kX^pY^q.
\]
\end{corollary}

\begin{proof}
Substitute \(X=cx\), \(Y=cy\), and \(s=t/c^2\) into
Theorem~\ref{thm:two-layer-kernel} and simplify.  The last identity follows
by expanding
\[
 \frac1{(1-X)(1-Y)-sXY}
 =\sum_{k\geq0}\frac{s^kX^kY^k}
 {(1-X)^{k+1}(1-Y)^{k+1}}.
\]
For the degree bound, also expand
\[
 \frac1{1-XY+sX^2Y^2}
 =\sum_{k\geq0}\frac{(-s)^kX^{2k}Y^{2k}}{(1-XY)^{k+1}}.
\]
In both inverse factors, each power of $s$ uses at least as many powers
of $X$ and of $Y$. Every monomial in $XY(1-XY)N(X,Y;s)$ has the
same property, and the remaining geometric factors use no $s$.
Thus $[s^jX^aY^b]\mathcal G=0$ whenever $j>\min(a,b)$.
\end{proof}

Coefficient extraction gives, for \(1\leq a\leq b\),
\begin{align*}
 \gamma_1(H_{(a,b)})&=a(b-1),\\
 \gamma_2(H_{(a,b)})
 &=\binom a2\binom{b-1}{2}
   +\binom a3+\binom{a-1}{2},\\
 [s^a]\Gamma_{a,b}(s)&=\binom{b-1}{a}.
\end{align*}
The first formula also follows at once from
Proposition~\ref{prop:first-gamma}.

\subsection{The first interlacing range}

\begin{theorem}
\label{thm:min-two-interlacing}
If \(\beta=(a,b)\) and \(\min(a,b)\leq2\), then
\[
 H_{\beta^-}(t)\preceq H_\beta(t).
\]
In particular, \(H_\beta(t)\) is real-rooted, and all its zeros are negative.
\end{theorem}

\begin{proof}
Coefficient extraction from Theorem~\ref{thm:two-layer-kernel} gives, for
\(b\geq2\),
\begin{equation}
 H_{(1,b)}(t)=c^b+(b-1)tc^{b-2}
 =c^{b-2}\bigl(t^2+(b+1)t+1\bigr).
 \label{eq:H-one-b}
\end{equation}
The initial case is $H_{(1,1)}=1+t$. The two roots different from $-1$ are negative and reciprocal, and they move
monotonically away from \(-1\) as \(b\) increases.  This proves
\(H_{(1,b-1)}\preceq H_{(1,b)}\).

For \(b\geq3\), the same extraction gives
\begin{equation}
 H_{(2,b)}(t)=c^{b+1}+2(b-1)tc^{b-1}
 +\binom{b-1}{2}t^2c^{b-3},
 \label{eq:H-two-b}
\end{equation}
while \(H_{(2,2)}(t)=1+5t+5t^2+t^3\).  Put \(x=b-1\).  After removing
the factor \(c^{x-2}\), divide the reciprocal quartic in
\eqref{eq:H-two-b} by \(t^2\) and set \(z=c^2/t\).  The result is
\[
 z^2+2xz+\binom x2,
\]
whose zeros are
\[
 z_x^\pm=-x\pm\sqrt{\frac{x(x+1)}2}<0.
\]
For \(x\geq3\), direct comparison gives
\[
 z_x^-<z_{x-1}^-<z_x^+<z_{x-1}^+<0.
\]
Every negative \(z\) gives a reciprocal pair of negative \(t\)-zeros, and
the additional zero \(-1\) at each growth step completes the alternation.
The initial step follows directly from the displayed polynomials.  Therefore
\(H_{(2,b-1)}\preceq H_{(2,b)}\).

It remains to compare the boundary steps after reversal. The cases
$a=1,2$ follow from $H_{(1)}=1$, $H_{(2)}=H_{(1,1)}=c$,
$H_{(1,2)}=c^2+t$, and $H_{(2,2)}=c^3+2tc$.
For $a\geq3$, the non-unit pair
of \(H_{(1,a)}\) in \eqref{eq:H-one-b} corresponds to \(z=-(a-1)\), which
lies strictly between \(z_{a-1}^-\) and \(z_{a-1}^+\).  Since
\(H_{(a)}(t)=c^{a-1}\), this proves
\[
 H_{(a)}\preceq H_{(a,1)}
 \quad\text{and}\quad
 H_{(a,1)}\preceq H_{(a,2)}.
\]
Reversal is an isomorphism of composition intervals, so these comparisons
cover every case with \(\min(a,b)\leq2\).
\end{proof}

The natural next step is to control the diagonal factor
\(1-XY+sX^2Y^2\) in Corollary~\ref{cor:two-layer-gamma-kernel}.  The other
factor is already a shifted Jacobi family.  Thus proving real-rootedness and
parameter interlacing for all \(\Gamma_{a,b}\) is now a concrete boundary
perturbation problem rather than an unstructured computation over
compositions.



\subsection{The first three rows and columns}

The first two coefficient formulas following
Corollary~\ref{cor:two-layer-gamma-kernel}, together with the formula for
the leading coefficient, determine the complete gamma polynomials in the
first three rows.  For \(x\geq2\), put
\begin{align}
 F_x(s)&=1+2xs+\binom{x}{2}s^2,\label{eq:F-x}\\
 G_x(s)&=1+3xs+\left(3\binom{x}{2}+2\right)s^2
              +\binom{x}{3}s^3.\label{eq:G-x}
\end{align}
Thus \(F_x=\Gamma_{2,x+1}\) and \(G_x=\Gamma_{3,x+1}\).

\begin{theorem}
\label{thm:min-three-interlacing}
If \(\beta=(a,b)\) and \(\min(a,b)\leq3\), then
\[
 H_{\beta^-}(t)\preceq H_\beta(t).
\]
In particular, \(H_\beta(t)\) has only negative real zeros.
\end{theorem}

\begin{proof}
The cases with \(\min(a,b)\leq2\) were proved above.  We first treat growth
along the third row.  From \eqref{eq:G-x}, for \(y\geq2\),
\begin{equation}
 G_{y+1}(s)-G_y(s)=sA_y(s),
 \qquad
 A_y(s)=3+3ys+\binom{y}{2}s^2.
 \label{eq:G-increment}
\end{equation}
The discriminant of \(A_y\) is \(3y(y+2)>0\), so its two zeros
\(\rho_1<\rho_2\) are negative.  At a zero of \(A_y\), polynomial division
gives
\begin{equation}
 G_y(s)=-\frac{y^2+3y+14}{y-1}s
        -\frac{2(y^2+2y+6)}{y(y-1)}.
 \label{eq:G-remainder}
\end{equation}
The zero of the right-hand side is
\[
 r_y=-\frac{2(y^2+2y+6)}{y(y^2+3y+14)}.
\]
A direct substitution gives
\[
 A_y(r_y)=
 -\frac{(y-1)(y^4+7y^3+28y^2+36y-72)}
 {y(y^2+3y+14)^2}<0.
\]
Hence \(\rho_1<r_y<\rho_2\).  It follows from
\eqref{eq:G-remainder} that \(G_y(\rho_1)>0>G_y(\rho_2)\), while
\(G_y(0)>0\).  For \(y\geq3\), the positive leading coefficient of \(G_y\)
also gives \(G_y(s)<0\) for all sufficiently negative \(s\).  These signs
place one zero of \(G_y\) in each of
\((-\infty,\rho_1)\), \((\rho_1,\rho_2)\), and \((\rho_2,0)\).
When \(y=2\), the last two intervals contain the two zeros of \(G_2\).
Consequently,
\[
 G_y\prec sA_y.
\]
The two polynomials on the right of \eqref{eq:G-increment} therefore have
interlacing zeros.  The convex-cone property of proper position now yields
\begin{equation}
 G_y\prec G_{y+1}\prec sA_y.
 \label{eq:G-consecutive-interlacing}
\end{equation}

It remains to compare the second and third rows, which is needed when the
last part grows from \(2\) to \(3\).  Equations \eqref{eq:F-x} and
\eqref{eq:G-x} give
\begin{equation}
 G_x(s)-F_x(s)=sB_x(s),
 \qquad
 B_x(s)=x+(x(x-1)+2)s+\binom{x}{3}s^2.
 \label{eq:FG-increment}
\end{equation}
For \(x=2\), the zeros of \(F_2\) are \(-2\pm\sqrt3\), and these alternate
with the zeros \(-1/2,0\) of \(sB_2\).  Suppose \(x\geq3\).  The
discriminant of \(B_x\) is
\[
 \frac{x^4+11x^2-12x+12}{3}>0,
\]
so \(B_x\) has two negative zeros \(\eta_1<\eta_2\).  At a zero of \(B_x\),
polynomial division gives
\[
 F_x(s)=-\frac{x^2+x+6}{x-2}s-\frac{2(x+1)}{x-2}.
\]
The zero of this linear remainder is
\[
 q_x=-\frac{2(x+1)}{x^2+x+6},
\]
and
\[
 B_x(q_x)=
 -\frac{(x-2)(x^4+4x^3+17x^2+14x-36)}
 {3(x^2+x+6)^2}<0.
\]
Thus \(\eta_1<q_x<\eta_2\), and the same sign argument shows
\(F_x\prec sB_x\).  Equation \eqref{eq:FG-increment} and the convex-cone
property give
\begin{equation}
 F_x\prec G_x.
 \label{eq:FG-interlacing}
\end{equation}

Relations \eqref{eq:G-consecutive-interlacing} and
\eqref{eq:FG-interlacing}, together with reversal symmetry of composition
intervals, cover every canonical step with minimum part equal to \(3\).
The gamma--interlacing correspondence for palindromic polynomials converts
these strict gamma relations into \(H_{\beta^-}\preceq H_\beta\).  The weak
symbol allows the forced common zeros at \(-1\).
\end{proof}



\subsection{The fourth row and column}

We next extend the finite-width theorem by one full row and column.  Besides
giving a larger exact range, the proof illustrates that the resulting
quartics need not be solved: two short pseudo-division chains reduce all
root comparisons to signs of quadratics and linear polynomials.

We use the following elementary sign certificate twice.

\begin{lemma}[Cubic--quadratic sign certificate]
\label{lem:cubic-quadratic-certificate}
Let \(P\) be a cubic with positive constant and leading coefficients, and
let \(R=r_0+r_1s+r_2s^2\), where \(r_0,r_1,r_2>0\) and \(R\) has two
distinct real zeros.  Suppose that, at every zero \(\rho\) of \(R\),
\[
 r_2^2P(\rho)=cL(\rho),
 \qquad L(s)=\ell_0+\ell_1s,
\]
where \(c>0\), \(\ell_0,\ell_1<0\), and
\(R(-\ell_0/\ell_1)<0\).  Then \(P\) has three distinct negative zeros and
\(R\) strictly interlaces \(P\).
\end{lemma}

\begin{proof}
The two zeros \(\rho_1<\rho_2\) of \(R\) are negative.  The zero
\(q=-\ell_0/\ell_1\) of \(L\) is also negative, and the last hypothesis says
\(\rho_1<q<\rho_2\).  Since \(\ell_1<0\), we have
\(P(\rho_1)>0>P(\rho_2)\).  Together with
\(P(-\infty)<0<P(0)\), these signs put one zero of \(P\) in each of
\((-\infty,\rho_1)\), \((\rho_1,\rho_2)\), and \((\rho_2,0)\).
\end{proof}

For \(x\geq3\), put
\begin{align}
 K_x(s)={}&1+4xs+\left(6\binom{x}{2}+7\right)s^2
 +\left(4\binom{x}{3}+2x-3\right)s^3
 +\binom{x}{4}s^4.
 \label{eq:K-x}
\end{align}
Thus \(K_x=\Gamma_{4,x+1}\), including the diagonal case \(x=3\).  Its
consecutive difference is
\begin{equation}
 K_{x+1}-K_x=sA_x,
 \qquad
 A_x=4+6xs+\left(4\binom{x}{2}+2\right)s^2
       +\binom{x}{3}s^3.
 \label{eq:K-increment}
\end{equation}
The discriminant of \(A_x\) is \(4D_A(x)\), where
\[
 D_A(x)=x^6+6x^5+21x^4+8x^3-96x^2+96x-32.
\]
After the substitution \(x=y+3\), every coefficient of \(D_A\) is positive.
Thus \(A_x\) has three distinct negative zeros.

Define
\begin{align*}
 R_x(s)&=r_0+r_1s+r_2s^2,\\
 r_0&=\frac{x^2+2x+12}{12},
 &r_1&=\frac{x(x^2+3x+20)}{12},
 &r_2&=\frac{x^4+3x^3+28x^2-4x+36}{72}.
\end{align*}
Its discriminant is
\[
 \frac{x^6+8x^5+55x^4+184x^3+472x^2-48x-864}{432}>0
 \qquad(x\geq3);
\]
the numerator again has only positive coefficients after \(x=y+3\).
Hence \(R_x\) has two distinct negative zeros.  Two steps of
pseudo-division give the following exact sign data.  At a zero \(\alpha\)
of \(A_x\),
\begin{equation}
 \binom{x}{3}^{\!2}K_x(\alpha)
 =-x(x-1)^2(x-2)R_x(\alpha),
 \label{eq:KA-pseudoremainder}
\end{equation}
whereas at a zero \(\rho\) of \(R_x\),
\begin{equation}
 r_2^2A_x(\rho)=x(x-2)L_x(\rho),
 \qquad
 L_x(s)=-\frac{U(x)}{2592}-\frac{V(x)}{5184}s,
 \label{eq:AR-pseudoremainder}
\end{equation}
with
\begin{align*}
 U(x)&=x^6+8x^5+63x^4+172x^3+348x^2+904x+936,\\
 V(x)&=x^7+10x^6+89x^5+336x^4+1056x^3
       +2036x^2+3480x-864.
\end{align*}
Both \(U(x)\) and \(V(x)\) are positive for \(x\geq3\).  If
\(q_x=-2U(x)/V(x)\) is the zero of \(L_x\), direct substitution and
factorization give
\begin{equation}
 R_x(q_x)\left(\frac{-V(x)}{5184}\right)^{\!2}
 =-\frac{W(x)^2P(x)}{967458816},
 \label{eq:row-four-location}
\end{equation}
where
\begin{align*}
 W(x)&=x^4+3x^3+28x^2-4x+36,\\
 P(x)&=x^8+16x^7+153x^6+842x^5+2756x^4+4824x^3
       +1904x^2-15328x-69408.
\end{align*}
Every coefficient of \(P(y+3)\) is positive, so the right-hand side of
\eqref{eq:row-four-location} is negative.  Lemma
\ref{lem:cubic-quadratic-certificate} now shows that \(R_x\) strictly
interlaces \(A_x\).

Let the zeros of \(A_x\) be
\(\alpha_1<\alpha_2<\alpha_3\).  The signs of \(R_x\) there are
\((+,-,+)\), and \eqref{eq:KA-pseudoremainder} therefore gives
\[
 K_x(\alpha_1)<0<K_x(\alpha_2),
 \qquad K_x(\alpha_3)<0.
\]
Together with the signs at \(-\infty\) and \(0\), this proves
\begin{equation}
 K_x\prec sA_x
 \qquad(x\geq3).
 \label{eq:K-A-proper-position}
\end{equation}
For \(x=3\), the three zeros of \(K_x\) occupy the last three intervals cut
out by the four zeros of \(sA_x\); for \(x\geq4\), the four zeros of \(K_x\)
and \(sA_x\) alternate.  Thus \eqref{eq:K-increment} and the convex-cone
property of proper position yield
\begin{equation}
 K_x\prec K_{x+1}.
 \label{eq:K-consecutive-interlacing}
\end{equation}

It remains to pass from the third row to the fourth.  With \(G_x\) as in
\eqref{eq:G-x}, we have
\begin{equation}
 K_x-G_x=sB_x,
 \quad
 B_x=x+\left(3\binom{x}{2}+5\right)s
 +\left(3\binom{x}{3}+2x-3\right)s^2
 +\binom{x}{4}s^3.
 \label{eq:GK-increment}
\end{equation}
When \(x=3\),
\[
 G_3=1+9s+11s^2+s^3,
 \qquad B_3=3+14s+6s^2.
\]
At a zero of \(B_3\), the remainder of \(G_3\) is
\(-(60+211s)/18\).  Its zero is \(-60/211\), and
\(B_3(-60/211)=-22077/44521<0\).  The resulting signs show directly that
\(G_3\prec sB_3\).

Suppose \(x\geq4\).  Define
\begin{align*}
 \widehat R_x(s)&=\widehat r_0+\widehat r_1s+\widehat r_2s^2,\\
 \widehat r_0&=\frac{x+1}{8},
 &\widehat r_1&=\frac{3x^2+3x+20}{24},
 &\widehat r_2&=\frac{x^3+11x-12}{48}.
\end{align*}
Its discriminant is
\[
 \frac{3x^4+12x^3+63x^2+126x+472}{576}>0.
\]
At a zero \(\delta\) of \(B_x\), pseudo-division gives
\begin{equation}
 \binom{x}{4}G_x(\delta)
 =-x(x-1)(x-2)\widehat R_x(\delta),
 \label{eq:GB-pseudoremainder}
\end{equation}
and at a zero \(\rho\) of \(\widehat R_x\), it gives
\begin{equation}
 \widehat r_2^{\,2}B_x(\rho)
 =(x-1)(x-3)\widehat L_x(\rho),
 \qquad
 \widehat L_x(s)=-\frac{\widehat U(x)}{4608}
                  -\frac{\widehat V(x)}{27648}s,
 \label{eq:BR-pseudoremainder}
\end{equation}
where
\begin{align*}
 \widehat U(x)&=x^5+4x^4+33x^3+34x^2+4x+144,\\
 \widehat V(x)&=3x^6+15x^5+147x^4+261x^3
                  +934x^2+1480x+2880.
\end{align*}
The zero \(\widehat q_x=-6\widehat U(x)/\widehat V(x)\) of
\(\widehat L_x\) satisfies
\begin{equation}
 \widehat R_x(\widehat q_x)
 \left(\frac{-\widehat V(x)}{27648}\right)^{\!2}
 =-\frac{(x-1)(x^2+x+12)^2Q(x)}{6115295232},
 \label{eq:boundary-four-location}
\end{equation}
where
\[
 Q(x)=3x^8+30x^7+264x^6+1110x^5+3441x^4+5916x^3
      +6196x^2-10592x-67968.
\]
Every coefficient of \(Q(y+3)\) is positive.  Lemma
\ref{lem:cubic-quadratic-certificate} therefore shows that
\(\widehat R_x\) strictly interlaces \(B_x\).  Equation
\eqref{eq:GB-pseudoremainder} gives the sign pattern \((-,+,-)\) for \(G_x\)
at the three zeros of \(B_x\).  Hence \(G_x\prec sB_x\), and
\eqref{eq:GK-increment} gives
\begin{equation}
 G_x\prec K_x
 \qquad(x\geq3).
 \label{eq:GK-interlacing}
\end{equation}

\begin{theorem}
\label{thm:min-four-interlacing}
If \(\beta=(a,b)\) and \(\min(a,b)\leq4\), then
\[
 H_{\beta^-}(t)\preceq H_\beta(t).
\]
In particular, \(H_\beta(t)\) has only negative real zeros.
\end{theorem}

\begin{proof}
Theorem~\ref{thm:min-three-interlacing} handles the first three rows and
columns.  Equations \eqref{eq:K-consecutive-interlacing} and
\eqref{eq:GK-interlacing} handle growth along the fourth row and the boundary
step from the third row to the fourth.  Reversal symmetry covers the fourth
column.  The gamma--interlacing correspondence converts these strict gamma
relations into weak interlacing for the characteristic Chow polynomials; the
weak symbol allows the forced common zeros at \(-1\).
\end{proof}



\subsection{Gamma differences and the diagonal boundary}

Continue to write
\[
 \mathcal G(X,Y;s)=\sum_{a,b\geq1}\Gamma_{a,b}(s)X^aY^b
\]
and put
\[
 Q(X,Y;s)=1-X-Y+(1-s)XY,
 \qquad F(Z;s)=1-Z+sZ^2.
\]
The rational gamma kernel admits a useful boundary decomposition.

\begin{proposition}[Delannoy propagation away from the diagonal]
\label{prop:gamma-delannoy-propagation}
For \(a,b\geq2\) with \(a\neq b\),
\begin{equation}
 \Gamma_{a,b}
 =\Gamma_{a-1,b}+\Gamma_{a,b-1}
 +(s-1)\Gamma_{a-1,b-1}.
 \label{eq:gamma-delannoy-recurrence}
\end{equation}
Let \(f_0(s)=f_1(s)=1\) and
\[
 f_m(s)=f_{m-1}(s)-s f_{m-2}(s),
 \qquad
 \sum_{m\geq0}f_m(s)Z^m=\frac1{F(Z;s)}.
\]
On the diagonal, the defect in \eqref{eq:gamma-delannoy-recurrence} is
\begin{equation}
 \Gamma_{a,a}-2\Gamma_{a-1,a}
 -(s-1)\Gamma_{a-1,a-1}=-s f_{a-1}(s)
 \qquad(a\geq2).
 \label{eq:gamma-diagonal-defect}
\end{equation}
\end{proposition}

\begin{proof}
Multiplying the rational expression for \(\mathcal G\) by \(Q\), and
simplifying, gives
\begin{equation}
 Q\mathcal G
 =XY+\frac{sXY^2}{1-Y}+\frac{sX^2Y}{1-X}
 -sXY\left(\frac1{F(XY;s)}-1\right).
 \label{eq:QG-boundary-support}
\end{equation}
The right-hand side is supported on the first row, the first column, and the
main diagonal.  The coefficient of \(X^aY^b\) on the left is precisely the
left-hand side of \eqref{eq:gamma-delannoy-recurrence}.  Hence it vanishes
for \(a,b\geq2\), \(a\neq b\).  The coefficient of \(X^aY^a\), for
\(a\geq2\), is \(-s f_{a-1}(s)\), proving
\eqref{eq:gamma-diagonal-defect}.
\end{proof}

For canonical growth in the second part, use the boundary convention
$H_{(a,0)}:=H_{(a)}=(1+t)^{a-1}$ and define
\[
 \Delta_{a,b}(s)=\Gamma_{a,b}(s)-\Gamma_{a,b-1}(s),
 \qquad \Gamma_{a,0}(s)=1.
\]

\begin{proposition}[Gamma-difference kernel]
\label{prop:gamma-difference-kernel}
We have
\begin{equation}
 \sum_{a,b\geq1}\Delta_{a,b}(s)X^aY^b
 =\frac{sXY\bigl(X+Y-2XY+sX^2Y^2\bigr)}
 {(1-X)F(XY;s)Q(X,Y;s)}.
 \label{eq:gamma-difference-kernel}
\end{equation}
If \(b>a\geq2\), then
\begin{align}
 \Delta_{a,b}(s)
 &=\Delta_{a-1,b}(s)+s\Gamma_{a-1,b-1}(s),
 \label{eq:gamma-difference-recursion}\\
 &=s\left(1+\sum_{r=1}^{a-1}\Gamma_{r,b-1}(s)\right).
 \label{eq:gamma-difference-prefix-sum}
\end{align}
The prefix-sum formula also holds for $a=1$, as $\Delta_{1,b}=s$ for $b>1$.
\end{proposition}

\begin{proof}
By definition,
\[
 \sum_{a,b\geq1}\Delta_{a,b}X^aY^b
 =(1-Y)\mathcal G(X,Y;s)-\frac{XY}{1-X}.
\]
Substitution of the rational expression for \(\mathcal G\) gives
\eqref{eq:gamma-difference-kernel}.  Rearranging
\eqref{eq:gamma-delannoy-recurrence} at $(a,b)$ gives
\eqref{eq:gamma-difference-recursion}.  Since
\(\Delta_{1,b}=s\) for \(b>1\), iteration in \(a\) yields
\eqref{eq:gamma-difference-prefix-sum}.
\end{proof}

The diagonal boundary itself has a compact Legendre-type generating
function.

\begin{theorem}[Diagonal gamma generating function]
\label{thm:diagonal-gamma-generating-function}
With \(\Gamma_{0,0}(s)=1\),
\begin{equation}
 \sum_{n\geq0}\Gamma_{n,n}(s)Z^n
 =\frac{(1-Z)(1-sZ)}{1-Z+sZ^2}
 \frac1{\sqrt{1-2(1+s)Z+(1-s)^2Z^2}}.
 \label{eq:diagonal-gamma-generating-function}
\end{equation}
\end{theorem}

\begin{proof}
Take the diagonal of the rational function \(\mathcal G\).  The identity
\[
 N(X,Y;s)=\bigl((1-s)-sXY\bigr)Q(X,Y;s)+s(1-sXY)
\]
reduces the calculation to
\[
 J(Z;s)=\operatorname{Diag}
 \frac1{(1-X)(1-Y)Q(X,Y;s)}.
\]
Expanding the last function gives
\[
 [Z^n]J(Z;s)
 =\sum_{k=0}^n\binom{n+1}{k+1}^2s^k
 =\frac{D_{n+1}(s)-1}{s},
\]
where
\[
 D_n(s)=\sum_k\binom nk^2s^k,
 \qquad
 \sum_{n\geq0}D_n(s)Z^n
 =\frac1{\sqrt{1-2(1+s)Z+(1-s)^2Z^2}}.
\]
Substitution and cancellation leave exactly
\eqref{eq:diagonal-gamma-generating-function}.
\end{proof}

Formula \eqref{eq:diagonal-gamma-generating-function} identifies the
remaining boundary as a rational perturbation of the central
Delannoy--Legendre family.  It suggests the following strengthening of
two-layer real-rootedness.

\begin{conjecture}[Gamma-difference Sturm property]
\label{conj:gamma-difference-sturm}
For all \(a,b\geq1\) with $(a,b)\neq(1,1)$, the nonzero polynomial
\(\Delta_{a,b}(s)\) has nonnegative coefficients, only nonpositive real zeros, and
\[
 \Gamma_{a,b-1}(s)\preceq\Delta_{a,b}(s).
\]
This gamma statement implies that the canonical Chow increment
\[
 C_{a,b}(t)=H_{(a,b)}(t)-H_{(a,b-1)}(t)
\]
is real-rooted and satisfies
\[
 H_{(a,b-1)}(t)\preceq C_{a,b}(t)
 \preceq H_{(a,b)}(t).
\]
\end{conjecture}

At $(a,b)=(1,1)$ the difference $\Delta_{1,1}$ is zero; the Chow assertion
holds directly since $H_{(1,0)}=1$, $C_{1,1}=t$, and $H_{(1,1)}=1+t$.
For all other pairs, the implication is Lemma~\ref{lem:gamma-increment-transfer},
with $p=\Gamma_{a,b-1}$, $g=\Delta_{a,b}$, and $n=a+b-2$.
The assumed nonnegativity of all differences, starting from $\Gamma_{a,0}=1$,
gives nonnegative coefficients for $p$, so its real zeros are negative.
The degree bounds follow from the ambient palindromic degrees, and
$g(0)=0$ because both gamma polynomials have constant term $1$.
Indeed, its decomposition is
\[
 C_{a,b}(t)
 =(1+t)^{a+b-1}\Delta_{a,b}\!\left(\frac{t}{(1+t)^2}\right)
 +t(1+t)^{a+b-2}\Gamma_{a,b-1}\!\left(\frac{t}{(1+t)^2}\right).
\]
Exact coefficient generation followed by numerical root checks finds no
failure of Conjecture~\ref{conj:gamma-difference-sturm} for
\(1\leq a,b\leq16\).  In addition, the diagonal sequence
\(\Gamma_{n,n}\) is a Sturm sequence throughout the tested range
\(n\leq24\). The diagonal check was also confirmed by exact rational root
isolation. These finite audits motivate a proof strategy but are not used as
a proof.


\subsection{Exact fixed-width transfer and Laguerre asymptotics}

The off-diagonal recurrence can be solved exactly, rather than only
iterated coefficient by coefficient.  Put \(\Gamma_{0,b}(s)=1\), and write
\(c_{a,k}(b)=[s^k]\Gamma_{a,b}(s)\).

\begin{proposition}[Exact fixed-width transfer]
\label{prop:exact-fixed-width-transfer}
Fix \(a\geq1\).  For \(b\geq a\),
\begin{equation}
 \Gamma_{a,b}(s)=\Gamma_{a,a}(s)
 +\sum_{j=0}^{a-1}\sum_{r=1}^{a-j}
 \binom{b-a}{r}\binom{a-j-1}{r-1}s^r\Gamma_{j,a}(s).
 \label{eq:exact-fixed-width-transfer}
\end{equation}
Consequently, for \(0\leq k\leq a\), \(c_{a,k}(b)\) is a polynomial in
\(b\) of degree at most \(k\) throughout the full ray \(b\geq a\).
More explicitly,
\begin{equation}
 c_{a,k}(b)=c_{a,k}(a)
 +\sum_{j=0}^{a-1}\sum_{r=1}^{\min(a-j,k)}
 \binom{b-a}{r}\binom{a-j-1}{r-1}c_{j,k-r}(a).
 \label{eq:exact-fixed-width-coefficients}
\end{equation}
\end{proposition}

\begin{proof}
For \(b>a\), equation \eqref{eq:gamma-difference-prefix-sum} gives, for
every \(0\leq i\leq a\),
\[
 \Gamma_{i,b}=\Gamma_{i,b-1}
 +s\sum_{j=0}^{i-1}\Gamma_{j,b-1}.
\]
Let \(N_a\) be the strictly lower triangular matrix whose entries below the
diagonal are all \(1\).  In vector form this recurrence is
\[
 (\Gamma_{0,b},\ldots,\Gamma_{a,b})^{\mathsf T}
 =(I+sN_a)(\Gamma_{0,b-1},\ldots,\Gamma_{a,b-1})^{\mathsf T}.
\]
After \(b-a\) steps, expand \((I+sN_a)^{b-a}\).  For \(i>j\) and
\(r\geq1\),
\[
 (N_a^r)_{i,j}=\binom{i-j-1}{r-1},
\]
because the entry counts strictly increasing chains of \(r-1\) intermediate
indices between \(j\) and \(i\).  Taking the last component proves
\eqref{eq:exact-fixed-width-transfer}; coefficient extraction gives
\eqref{eq:exact-fixed-width-coefficients}.
\end{proof}

The exact transfer also gives the first correction to the Laguerre limit.

We first record the elementary analytic step that converts a scaled limit
into consecutive interlacing. It will also be used for the sandwich family.

\begin{lemma}[Eventual interlacing from analytic scaling]
\label{lem:analytic-scaling-interlacing}
Let $p_m(z)=\sum_{j=0}^d a_j(m)z^j$, where each $a_j$ is a real
polynomial of degree at most $j$. Suppose that the coefficientwise limit
$P(z)=\lim_{m\to\infty}p_m(z/m)$ has degree $d$ and $d$ simple
negative zeros $\zeta_1<\cdots<\zeta_d$. Then, for all sufficiently
large real $m$, $p_m$ has degree $d$ and simple negative zeros
$\rho_1(m)<\cdots<\rho_d(m)$, with
\begin{align*}
 \rho_i(m)&=\zeta_i/m+O(m^{-2}),&
 \rho_i'(m)&=-\zeta_i/m^2+O(m^{-3})>0.
\end{align*}
If $d\geq1$, then, for all sufficiently large $m$,
\[
 \rho_i(m-1)<\rho_i(m)<\rho_{i+1}(m-1)\quad(1\leq i<d),
 \qquad \rho_d(m-1)<\rho_d(m)<0.
\]
For $d=0$, $p_m$ is a nonzero constant and the root-order assertions
are vacuous.
\end{lemma}

\begin{proof}
Set $F(\varepsilon,z)=\sum_j\varepsilon^j a_j(1/\varepsilon)z^j$.
The degree bounds make $F$ a polynomial in $\varepsilon$ and $z$,
with $F(0,z)=P(z)$. The real analytic implicit-function theorem at
each simple zero gives real branches $z_i(\varepsilon)$ with
$z_i(0)=\zeta_i$. For small real $\varepsilon$, these branches are
distinct and negative; the nonzero limiting leading coefficient
ensures that they account for the full degree.

Now $\rho_i(m)=z_i(1/m)/m$. Differentiating this exact analytic
expression gives
\[
 \rho_i'(m)=-\frac{z_i(1/m)}{m^2}
           -\frac{z_i'(1/m)}{m^3}
           =-\frac{\zeta_i}{m^2}+O(m^{-3}).
\]
Hence every root moves strictly to the right. Since there are finitely
many branches, one common threshold and constants $c,C>0$ give
$\rho_{i+1}(m-1)-\rho_i(m-1)\geq c/(m-1)$ and
$0<\rho_i(m)-\rho_i(m-1)\leq C/(m-1)^2$, the latter by the
mean-value theorem. Increasing the threshold so that $C/(m-1)<c$
proves the displayed cross-order. For $d=1$ only monotonicity and
negativity are needed; for $d=0$ the polynomial is a nonzero constant.
\end{proof}

\begin{theorem}[Eventual fixed-width interlacing]
\label{thm:eventual-fixed-width}
For every fixed \(a\geq1\), there is an integer \(B_a\geq a+2\) such that
\[
 H_{(a,b-1)}(t)\preceq H_{(a,b)}(t)
 \qquad (b\geq B_a).
\]
In particular, \(H_{(a,b)}(t)\) has only negative real zeros for all
sufficiently large \(b\).  More precisely, if
\[
 P_a(z)=L_a(-z)=\sum_{k=0}^{a}\binom{a}{k}\frac{z^k}{k!},
 \qquad \vartheta=z\frac{d}{dz},
\]
where \(L_a\) is the ordinary Laguerre polynomial, then
\begin{equation}
 \Gamma_{a,b}\left(\frac{z}{b}\right)
 =P_a(z)-\frac{1}{2b}(\vartheta^2+\vartheta)P_a(z)
 +O(b^{-2})
 \label{eq:second-order-laguerre-limit}
\end{equation}
coefficientwise as \(b\to\infty\).
\end{theorem}

\begin{proof}
Formula \eqref{eq:exact-fixed-width-coefficients} gives the two highest
powers of \(b\) explicitly.  The case \(k=1\) is the exact identity
\(c_{a,1}(b)=a(b-1)\).  For \(2\leq k\leq a\), put \(n=b-a\).  Only the
terms with \(r=k\) and \(r=k-1\) contribute to the two highest powers of
\(n\).  Since \(c_{j,1}(a)=j(a-1)\), the hockey-stick identity and
\[
 \sum_{j=1}^{a-k+1}j\binom{a-j-1}{k-2}=\binom{a}{k}
\]
give
\begin{equation}
 c_{a,k}(b)=\binom{a}{k}\frac{b^k}{k!}
 -\binom{a}{k}\frac{k+1}{2}\frac{b^{k-1}}{(k-1)!}
 +O(b^{k-2}).
 \label{eq:fixed-width-two-leading-coefficients}
\end{equation}
Together with the cases \(k=0,1\), substitution into
\(\Gamma_{a,b}(z/b)\) proves
\eqref{eq:second-order-laguerre-limit}.

Use the coefficient polynomials to extend $\Gamma_{a,b}$ to real $b$.
The degree bound in Corollary~\ref{cor:two-layer-gamma-kernel} shows
that there are no coefficients beyond $s^a$.
Let \(\zeta_1<\cdots<\zeta_a<0\) be the distinct zeros of \(P_a\).  The
Laguerre differential equation, in the present normalization, is
\[
 zP_a''(z)+(1+z)P_a'(z)-aP_a(z)=0.
\]
At a zero \(\zeta\) of \(P_a\), it follows that
\[
 -\frac12(\vartheta^2+\vartheta)P_a(\zeta)
 =-\frac{\zeta(1-\zeta)}2P_a'(\zeta).
\]
Because the coefficients in \eqref{eq:second-order-laguerre-limit} are
polynomial in \(1/b\), the implicit-function theorem gives analytic scaled
roots
\[
 z_j(b)=\zeta_j+\frac{\zeta_j(1-\zeta_j)}{2b}+O(b^{-2}).
\]
Thus the corresponding zeros of \(\Gamma_{a,b}\) satisfy
\begin{equation}
 \rho_j(b)=\frac{\zeta_j}{b}
 +\frac{\zeta_j(1-\zeta_j)}{2b^2}+O(b^{-3}).
 \label{eq:fixed-width-root-expansion}
\end{equation}
The coefficient bounds in
\eqref{eq:exact-fixed-width-coefficients} and the simple negative zeros
of $P_a$ verify every hypothesis of
Lemma~\ref{lem:analytic-scaling-interlacing}. It supplies the directed
order $\rho_j(b-1)<\rho_j(b)<\rho_{j+1}(b-1)$ for $j<a$, and the final inequality
$\rho_a(b-1)<\rho_a(b)<0$, with one threshold for all roots.
Choose $B_a\geq a+2$ so that also $2a\leq a+b-2$.
Applying Lemma~\ref{lem:gamma-chow-interlacing} with $n=a+b-2$
then proves the Chow interlacing, including its forced common zeros
at $-1$.
\end{proof}


\section{An arbitrary fixed number of layers}

The same-length convolution behind the two-layer kernel persists for every
fixed number of parts.  The only additional issue is that a shorter
composition may embed in several sets of positions.  A canonical rightmost
embedding resolves this ambiguity without inclusion--exclusion over
occurrences and leads to a finite rational operator.

For \(k\geq1\), let \(\mathbb P^k\) denote the positive compositions with
exactly \(k\) parts.  If
\(\alpha=(a_1,\ldots,a_k)\), write
\(\mathbf x^\alpha=x_1^{a_1}\cdots x_k^{a_k}\), and define
\[
 \mathcal R_k(\mathbf x;t)
 =\sum_{\alpha\in\mathbb P^k}R_\alpha(t)\mathbf x^\alpha,
 \qquad
 \mathcal H_k(\mathbf x;t)
 =\sum_{\beta\in\mathbb P^k}H_\beta(t)\mathbf x^\beta.
\]

We first record the rational operator that sums over all upper words.  Let
\[
 F(\mathbf z)=\sum_{\nu\in\mathbb P^j}f_\nu\mathbf z^\nu
\]
and, for \(0\leq j\leq k\), put
\begin{equation}
 \bigl(\mathscr U_{k,j}^{(q)}F\bigr)(\mathbf x)
 =\sum_{\nu\in\mathbb P^j}
   \sum_{\substack{\alpha\in\mathbb P^k\\\nu\leq\alpha}}
   f_\nu q^{|\alpha|-|\nu|}\mathbf x^\alpha.
 \label{eq:upper-word-operator}
\end{equation}
For \(j=0\), the inner word is empty, \(F=1\), and
\begin{equation}
 \mathscr U_{k,0}^{(q)}1
 =\prod_{p=1}^k\frac{qx_p}{1-qx_p}.
 \label{eq:empty-upper-word-operator}
\end{equation}

\begin{lemma}[Rightmost-embedding operator]
\label{lem:rightmost-embedding-operator}
If \(F\) is rational, then
\(\mathscr U_{k,j}^{(q)}F\) is rational.
\end{lemma}

\begin{proof}
Suppose \(j\geq1\), and let
\(I=(i_1<\cdots<i_j)\) be a possible rightmost embedding.  Set
\[
 G_0(I)=\{p:p<i_1\},\qquad
 G_r(I)=\{p:i_r<p<i_{r+1}\}\quad(1\leq r\leq j),
\]
where \(i_{j+1}=k+1\).  For subsets \(S_r\subseteq G_r(I)\), define
\[
 z_r(I,\mathbf S)
 =q^{|S_r|}x_{i_r}\prod_{p\in S_r}x_p.
\]
Then
\begin{align}
 \mathscr U_{k,j}^{(q)}F
 ={}&\sum_I
 \frac{\displaystyle\prod_{p\in G_0(I)}qx_p}
      {\displaystyle\prod_{p=1}^k(1-qx_p)}
 \sum_{\substack{S_r\subseteq G_r(I)\\1\leq r\leq j}}
 (-1)^{\sum_r|S_r|}
 \left(\prod_{r=1}^j
       \prod_{p\in G_r(I)\setminus S_r}qx_p\right)
 \notag\\[-2pt]
 &\hspace{42mm}\times
 F\bigl(z_1(I,\mathbf S),\ldots,z_j(I,\mathbf S)\bigr).
 \label{eq:rightmost-embedding-operator}
\end{align}
This is a finite sum of rational substitutions, so it remains only to
justify the formula.

For each pair \(\nu\leq\alpha\), choose the embedding greedily from the
right: \(i_j\) is the last position capable of holding \(\nu_j\), and the
earlier positions are chosen recursively.  This gives a unique \(I\).
For fixed \(I\), positions in \(G_0(I)\) are unrestricted, the chosen
position \(i_r\) is at least \(\nu_r\), and every position in \(G_r(I)\) is
strictly less than \(\nu_r\).  Their respective generating factors are
\[
 \frac{qx_p}{1-qx_p},\qquad
 \frac{x_{i_r}^{\nu_r}}{1-qx_{i_r}},\qquad
 \frac{qx_p-(qx_p)^{\nu_r}}{1-qx_p}.
\]
Expanding the numerator of the last factor over each gap and then summing
over \(\nu\) gives exactly
\eqref{eq:rightmost-embedding-operator}.
\end{proof}

The contribution of lower words having the same number of parts is governed
by a universal kernel.  Let \(e_m(\mathbf x)\) be the \(m\)-th elementary
symmetric polynomial and set
\[
 \mathcal K_k(\mathbf x;t)
 =\sum_{\mathbf d\in\mathbb N^k}
 [|\mathbf d|-1]_t\mathbf x^{\mathbf d}.
\]

\begin{proposition}[Universal fixed-layer kernel]
\label{prop:universal-fixed-layer-kernel}
For every \(k\geq1\),
\begin{equation}
 \mathcal K_k(\mathbf x;t)
 =\frac{
 t\left(\displaystyle\prod_{p=1}^k\frac1{1-x_p}-1\right)
 -\left(\displaystyle\prod_{p=1}^k\frac1{1-tx_p}-1\right)}
 {t(1-t)}.
 \label{eq:universal-K}
\end{equation}
Moreover,
\begin{equation}
 1-t\mathcal K_k(\mathbf x;t)
 =\frac{D_k(\mathbf x;t)}
 {\displaystyle\prod_{p=1}^k(1-x_p)(1-tx_p)},
 \qquad
 D_k(\mathbf x;t)
 =\sum_{m=0}^k(-1)^m[m+1]_t e_m(\mathbf x).
 \label{eq:universal-denominator}
\end{equation}
\end{proposition}

\begin{proof}
For a monomial of positive total degree \(d\geq1\), the coefficient on the right-hand
side of \eqref{eq:universal-K} is
\[
 \frac{t-t^d}{t(1-t)}=[d-1]_t.
\]
The constant coefficient is zero because the two products in
\eqref{eq:universal-K} each have their constant term removed; the case
$d=1$ also vanishes, as required.  Clearing the two product
denominators in \(1-t\mathcal K_k\) gives
\[
 D_k(\mathbf x;t)
 =\frac{\prod_p(1-x_p)-t\prod_p(1-tx_p)}{1-t}.
\]
Expansion in elementary symmetric polynomials proves the second expression
in \eqref{eq:universal-denominator}.
\end{proof}

\begin{theorem}[Fixed-layer rationality]
\label{thm:fixed-layer-rationality}
For every fixed \(k\geq1\), both
\(\mathcal R_k(\mathbf x;t)\) and
\(\mathcal H_k(\mathbf x;t)\) are rational functions.
More precisely, define
\begin{align}
 E_k(\mathbf x;t)
 &=
 \frac{\mathscr U_{k,0}^{(1)}1
       -t^{-1}\mathscr U_{k,0}^{(t)}1}{1-t},
 \label{eq:empty-contribution}\\
 A_k(\mathbf x;t)
 &=E_k(\mathbf x;t)
 +\sum_{j=1}^{k-1}
 \frac{
 t\mathscr U_{k,j}^{(1)}\mathcal R_j
 -\mathscr U_{k,j}^{(t)}\mathcal R_j}{1-t}.
 \label{eq:shorter-word-contribution}
\end{align}
Then
\begin{equation}
 \mathcal R_k(\mathbf x;t)
 =\frac{A_k(\mathbf x;t)}
 {1-t\mathcal K_k(\mathbf x;t)},
 \label{eq:fixed-layer-R-recursion}
\end{equation}
and
\begin{equation}
 \mathcal H_k(\mathbf x;t)
 =\mathscr U_{k,0}^{(1)}1
 +t\sum_{j=1}^k
 \mathscr U_{k,j}^{(1)}\mathcal R_j.
 \label{eq:fixed-layer-H-recursion}
\end{equation}
\end{theorem}

\begin{proof}
Apply the terminal recursion \eqref{eq:terminal-recursion} to every
\(\alpha\in\mathbb P^k\) and partition the lower compositions by their
number of parts.  The empty lower word contributes \(E_k\).  If
\(\nu\) has \(j<k\) parts and \(d=|\alpha|-|\nu|\), then \(d\geq1\) and
\[
 t[d-1]_t=\frac{t-t^d}{1-t}.
\]
Consequently, the total contribution from all \(j\)-part lower words is the
\(j\)-th summand in \eqref{eq:shorter-word-contribution}.

If \(\nu\) and \(\alpha\) both have \(k\) parts, then
\(\nu\leq\alpha\) is componentwise containment.  Writing
\(\alpha=\nu+\mathbf d\) turns this part of the terminal recursion into the
convolution
\[
 t\mathcal R_k(\mathbf x;t)\mathcal K_k(\mathbf x;t).
\]
This proves \eqref{eq:fixed-layer-R-recursion}.  Lemma
\ref{lem:rightmost-embedding-operator} and induction on \(k\) now prove that
\(\mathcal R_k\) is rational.  Finally, summing the terminal identity
\eqref{eq:terminal-sum} over all \(k\)-part upper words gives
\eqref{eq:fixed-layer-H-recursion}, and hence the rationality of
\(\mathcal H_k\).  Although
\eqref{eq:empty-contribution}--\eqref{eq:shorter-word-contribution} are
written in \(\mathbb Q(t,\mathbf x)\), their apparent factors
\(t^{-1}\) and \((1-t)^{-1}\) cancel coefficientwise by the displayed
\(t\)-integer identities.
\end{proof}

For three layers, the universal denominator is already quite compact:
\begin{equation}
 D_3(x,y,z;t)
 =1-[2]_t(x+y+z)+[3]_t(xy+xz+yz)-[4]_txyz.
 \label{eq:three-layer-denominator}
\end{equation}
Thus
\[
 \mathcal R_3(x,y,z;t)
 =\frac{(1-x)(1-y)(1-z)(1-tx)(1-ty)(1-tz)}
 {D_3(x,y,z;t)}\,A_3(x,y,z;t),
\]
where \(A_3\) is the finite rightmost-embedding expression obtained from
\(\mathcal R_1\) and \(\mathcal R_2\) in
\eqref{eq:shorter-word-contribution}.  The two-layer denominator from
\eqref{eq:kernel-cancellation} is the case
\(D_2=1-[2]_t(x+y)+[3]_txy\).  Hence the two-layer calculation is the first
nontrivial member of a fixed-layer rational hierarchy, rather than an
isolated special case.

\begin{remark}
The denominator $D_k$ is also the denominator of the augmented Chow
series associated with the Toeplitz matrix of complete homogeneous
symmetric functions: in \cite[Theorem~8.1]{BrandenVecchi2025TN}, substitute
$f(z)=\prod_{i=1}^k(1-x_i z)^{-1}$ to obtain
$\sum_{n\geq0}G_n(t)z^n=1/D_k(z\mathbf x;t)$.
The fixed-layer theorem adds the rightmost-embedding operator and the
recursive numerators for composition intervals. Real-rootedness under
nonnegative numerical specialization of $\mathbf x$ does not by itself
imply real-rootedness of each multivariate coefficient.
\end{remark}


\section{A rectangular family with unboundedly many layers}
\label{sec:singleton-prefix}

The fixed-layer rationality theorem does not merely repackage the two-layer
case.  The following family has an unbounded number of layers and admits a
closed form strong enough to settle real-rootedness and interlacing for every
parameter value.

\begin{theorem}[Singleton-prefix rectangles]
\label{thm:singleton-prefix-rectangles}
For $r\geq0$ and $b\geq1$, let
\(
  \beta_{r,b}=(1^r,b)
\)
and put $q=b-1$.  Define its gamma polynomial $\Delta_{r,q}(s)$ by
\[
  H_{(1^r,q+1)}(t)
  =(1+t)^{r+q}\Delta_{r,q}\!\left(\frac{t}{(1+t)^2}\right).
\]
Then the following statements hold.
\begin{enumerate}
\item The nonempty part of the lower interval is a rectangular grid:
\[
  P_{\beta_{r,b}}\setminus\{\varnothing\}
  \cong C_{r+1}\times C_b.
\]
\item The complete bivariate Chow series is
\begin{equation}
\label{eq:singleton-prefix-chow-series}
  \sum_{r,q\geq0}H_{(1^r,q+1)}(t)X^rY^q
  =\frac{1}{1-[2]_t(X+Y)+[3]_tXY}.
\end{equation}
Equivalently,
\begin{equation}
\label{eq:singleton-prefix-chow-formula}
  H_{(1^r,b)}(t)
  =\sum_{k=0}^{\min(r,b-1)}
    \binom{r}{k}\binom{b-1}{k}
    t^k(1+t)^{r+b-1-2k}.
\end{equation}
\item The gamma polynomial is the Delannoy--Jacobi polynomial
\begin{equation}
\label{eq:singleton-prefix-gamma}
  \Delta_{r,q}(s)
  =\sum_{k=0}^{\min(r,q)}\binom{r}{k}\binom{q}{k}s^k.
\end{equation}
Unless it is constant, it has only simple negative zeros.  Consequently,
every $H_{(1^r,b)}(t)$ is real-rooted.
\item For each fixed $r\geq1$, the consecutive polynomials
$\Delta_{r,q}$ and $\Delta_{r,q+1}$ strictly interlace for all $q\geq0$,
with the usual constant--linear convention at $q=0$.  Hence
$H_{(1^r,b)}(t)$ and $H_{(1^r,b+1)}(t)$ weakly interlace for every $b\geq1$;
equality can occur only at the forced root $t=-1$.
\end{enumerate}
\end{theorem}

\begin{proof}
Every nonempty composition below $(1^r,b)$ is of exactly one of the forms
\[
  1^{p+1}\quad(0\leq p\leq r),
  \qquad
  (1^p,j)\quad(0\leq p\leq r,\ 2\leq j\leq b).
\]
Indeed, only the final letter of the host can support a letter larger than
$1$, and such a letter must therefore be the final letter of the embedded
composition.  Sending $(p,1)$ to $1^{p+1}$ and $(p,j)$ to $(1^p,j)$ for
$j\geq2$ gives an order isomorphism from
$\{0,\ldots,r\}\times\{1,\ldots,b\}$ with componentwise order.  This proves
the first assertion.

We next compute the Chow series directly from the terminal-chain recurrence.
Let $R_{p,q}(t)$ be the terminal-chain polynomial attached to the grid vertex
$(p,q)$, whose rank in the interval is $p+q+1$, and set
\(
  \mathcal R=\sum_{p,q\geq0}R_{p,q}(t)X^pY^q.
\)
On a rectangular grid the recurrence becomes
\[
  \mathcal R=A+t\mathcal R\mathcal K_2,
  \qquad
  A=\sum_{p,q\geq0}[p+q]_tX^pY^q,
\]
where
\(
  \mathcal K_2=\sum_{(i,j)\in\mathbb N^2\setminus\{(0,0)\}}
  [i+j-1]_tX^iY^j.
\)
A routine summation gives
\begin{align*}
  A&=\frac{X+Y-(1+t)XY}
  {(1-X)(1-Y)(1-tX)(1-tY)},\\
  1-t\mathcal K_2
  &=\frac{1-[2]_t(X+Y)+[3]_tXY}
  {(1-X)(1-Y)(1-tX)(1-tY)}.
\end{align*}
Thus, with $D_2=1-[2]_t(X+Y)+[3]_tXY$,
\[
  \mathcal R=\frac{X+Y-(1+t)XY}{D_2}.
\]
Summing terminal contributions over every lower grid vertex now yields
\[
  \sum_{r,q\geq0}H_{(1^r,q+1)}(t)X^rY^q
  =\frac{1+t\mathcal R}{(1-X)(1-Y)}
  =\frac1{D_2},
\]
which is \eqref{eq:singleton-prefix-chow-series}.

Put $c=1+t$ and $s=t/c^2$.  Since
\(
  [3]_t=c^2-t=c^2(1-s),
\)
the gamma transform of \eqref{eq:singleton-prefix-chow-series} is
\[
  \sum_{r,q\geq0}\Delta_{r,q}(s)U^rV^q
  =\frac1{1-U-V+(1-s)UV}.
\]
Expanding the right-hand side as
\[
  \frac1{(1-U)(1-V)-sUV}
  =\sum_{k\geq0}
    \frac{s^kU^kV^k}{(1-U)^{k+1}(1-V)^{k+1}}
\]
gives \eqref{eq:singleton-prefix-gamma}, and its gamma expansion gives
\eqref{eq:singleton-prefix-chow-formula}.

It remains to locate the zeros.  Write
\(
  J_{a,c}(s)=\sum_k\binom ak\binom ck s^k.
\)
For $c\geq a$, the binomial expansion of the Jacobi polynomial gives
\begin{equation}
\label{eq:delannoy-jacobi}
  J_{a,c}(s)
  =(1-s)^aP_a^{(0,c-a)}\!\left(\frac{1+s}{1-s}\right).
\end{equation}
Indeed, after the substitution $x=(1+s)/(1-s)$, the two binomial factors
in the Jacobi expansion contribute $s^{a-m}/(1-s)^a$, and reindexing by
$k=a-m$ gives the defining sum for $J_{a,c}$.  The zeros of
$P_a^{(0,c-a)}$ are simple and lie in $(-1,1)$; under the increasing change
of variables $s=(x-1)/(x+1)$ they become simple negative zeros.  Symmetry in
$a$ and $c$ handles the remaining case.  The gamma-to-Chow transform then
proves real-rootedness of the Chow polynomials.

For completeness, interlacing follows from the standard Jacobi contiguous identity \cite[Chapter IV]{Szego1939}:
\begin{equation}
\label{eq:jacobi-contiguous}
  (2n+\lambda+1)P_n^{(0,\lambda)}
  =(n+\lambda+1)P_n^{(0,\lambda+1)}
    +nP_{n-1}^{(0,\lambda+1)}.
\end{equation}
The two polynomials on the right are consecutive members of the same Jacobi
family and have strictly interlacing zeros and positive leading coefficients.  Evaluating
\eqref{eq:jacobi-contiguous} at the zeros of either member shows that
$P_n^{(0,\lambda)}$ strictly interlaces the degree-$n$ member and is strictly
interlaced by the degree-$(n-1)$ member.

If $q<r$, use $J_{r,q}=J_{q,r}$ in \eqref{eq:delannoy-jacobi} and take
$n=q+1$ and $\lambda=r-q-1$ in \eqref{eq:jacobi-contiguous}.  If $q\geq r$,
take $n=r$ and $\lambda=q-r$.  These two cases prove strict interlacing of
$J_{r,q}$ and $J_{r,q+1}$ for every $r\geq1$ and $q\geq0$.  The monotone
change of variables in \eqref{eq:delannoy-jacobi} preserves the order of the
zeros. When the degrees agree, those of $J_{r,q+1}$ lie to the right of
the corresponding zeros of $J_{r,q}$; when the degree increases, its zeros
surround those of $J_{r,q}$. These orders follow from the signs in
\eqref{eq:jacobi-contiguous} and the positive leading coefficients.
Finally, Lemma~\ref{lem:gamma-chow-interlacing}
transfers this to weak interlacing of the Chow polynomials; the
only possible common zero is the power of $1+t$ forced by the gamma
expansion.
\end{proof}

\begin{remark}
The denominator in \eqref{eq:singleton-prefix-chow-series} is exactly the
universal two-coordinate kernel $D_2$.  Here it governs a family with
arbitrarily many layers because the interval collapses to a two-dimensional
grid.  Thus the two-coordinate kernel is structural rather than an artifact
of restricting the host to two layers.
\end{remark}


\section{Singleton-sandwich intervals}
\label{sec:singleton-sandwich}

The rectangular family remains tractable after singleton layers are allowed
on both sides of the unique nonsingleton layer.  The interval is no longer a
product of chains, but it splits into one chain and one three-dimensional
box.  This produces an explicit rational Chow series with the universal
three-coordinate denominator and, after the gamma transform, a manifestly
nonnegative series.

Put $c=1+t$ and define
\begin{equation}
\label{eq:sandwich-base-series}
\begin{aligned}
 \mathcal B(W,Z;t)
 :={}&\frac1{1-t}
 \left(
  \frac1{(1-W)(1-Z)}
  -\frac{t}{(1-tW)(1-tZ)}
 \right)\\
 &+\frac{tW}{(1-cW)(1-t)}
 \left(
  \frac1{(1-W)(1-Z)}
  -\frac1{(1-tW)(1-tZ)}
 \right).
\end{aligned}
\end{equation}
Further, let
\begin{equation}
\label{eq:sandwich-lifted-base}
 \mathcal A(X,Y,Z;t)
 :=\frac{X\mathcal B(X,Z;t)-Y\mathcal B(Y,Z;t)}{X-Y}
\end{equation}
and
\begin{equation}
\label{eq:sandwich-D3}
 D_3(X,Y,Z;t)
 :=1-[2]_t(X+Y+Z)+[3]_t(XY+XZ+YZ)-[4]_tXYZ.
\end{equation}
The divided difference in \eqref{eq:sandwich-lifted-base} is understood by
continuity when $X=Y$; its numerator is divisible by $X-Y$ as a formal
power series.

\begin{theorem}[Singleton-sandwich series]
\label{thm:singleton-sandwich-series}
The characteristic Chow polynomials of all hosts with one nonsingleton
layer satisfy
\begin{equation}
\label{eq:singleton-sandwich-series}
\begin{aligned}
 &\sum_{r,s,h\geq0}
 H_{(1^r,h+2,1^s)}(t)X^rY^sZ^h\\
 &\qquad=
 \frac1{(1-cX)(1-cY)(1-Z)}
 +\frac{t(1-tX)(1-tY)(1-tZ)}{D_3(X,Y,Z;t)}
  \mathcal A(X,Y,Z;t).
\end{aligned}
\end{equation}
In particular, this three-parameter family, whose number of layers is
unbounded, has a rational Chow generating function.

Define $\Delta_{r,s,h}(z)$ by
\[
 H_{(1^r,h+2,1^s)}(t)
 =c^{r+s+h+1}\Delta_{r,s,h}\!\left(\frac{t}{c^2}\right).
\]
Then the complete gamma series is
\begin{equation}
\label{eq:singleton-sandwich-gamma-series}
 \sum_{r,s,h\geq0}\Delta_{r,s,h}(z)U^rV^sW^h
 =\frac{1+\dfrac{zU}{1-U}+\dfrac{zV}{1-V}}
 {1-(U+V+W)+(1-z)(UV+UW+VW)-(1-2z)UVW}.
\end{equation}
Consequently, every $H_{(1^r,b,1^s)}(t)$ is gamma-positive.
\end{theorem}

\begin{proof}
Fix a host $(1^r,b,1^s)$ with $b\geq2$.  Every nonempty composition below it
is either $1^m$, with $1\leq m\leq r+s+1$, or has the form
\[
 (1^p,j,1^q),
 \qquad 0\leq p\leq r,\quad 2\leq j\leq b,\quad 0\leq q\leq s.
\]
Among the latter elements, generalized subword order is componentwise in
$(p,j,q)$.  Moreover,
\[
 1^m\leq(1^p,j,1^q)
 \quad\Longleftrightarrow\quad
 m\leq p+q+1.
\]
Thus the interval is the asserted chain-plus-box poset, including all of
its cross-relations.

The all-one interval is a chain, so
\[
 R_{(1)}(t)=0,
 \qquad
 R_{1^m}(t)=c^{m-2}\quad(m\geq2).
\]
Represent $(1^p,j,1^q)$ by the box coordinate $(p,h,q)$, where $h=j-2$,
and put $n=p+q$.  In its terminal recursion, the direct jump from the
bottom together with all penultimate elements on the all-one chain is
\begin{equation}
 B_{n,h}(t)
 =[n+h+1]_t
 +t\sum_{u=0}^{n-1}c^u[n+h-u-1]_t,
 \label{eq:sandwich-base-coefficient}
\end{equation}
where the sum is empty for $n=0$.  Summing the first term in
\eqref{eq:sandwich-base-coefficient} over $n,h\geq0$ gives the first line of
\eqref{eq:sandwich-base-series}.  In the second term, write
$n=u+d+1$; its series is
\[
 \frac{tW}{1-cW}
 \sum_{d,h\geq0}[d+h]_tW^dZ^h,
\]
which is the second line of \eqref{eq:sandwich-base-series}.

For fixed $n$, splitting the $n$ singleton positions as $n=p+q$ gives
\[
 \sum_{p+q=n}X^pY^q=\frac{X^{n+1}-Y^{n+1}}{X-Y}.
\]
Consequently, \eqref{eq:sandwich-lifted-base} is precisely the generating
function for the base terms $B_{p+q,h}$.

Let
\[
 \mathcal R_{\Box}(X,Y,Z;t)
 =\sum_{p,h,q\geq0}R_{(1^p,h+2,1^q)}(t)X^pY^qZ^h.
\]
All remaining penultimate elements lie in the componentwise box, so the
terminal recurrence is the three-coordinate convolution
\[
 \mathcal R_{\Box}
 =\mathcal A+t\mathcal R_{\Box}\mathcal K_3.
\]
By Proposition~\ref{prop:universal-fixed-layer-kernel},
\[
 \mathcal R_{\Box}
 =\mathcal A\,
 \frac{(1-X)(1-Y)(1-Z)(1-tX)(1-tY)(1-tZ)}
 {D_3(X,Y,Z;t)}.
\]

Finally, the empty element and all all-one elements contribute
$H_{1^{r+s+1}}(t)=c^{r+s}$.  Summing the box terminal contributions over
all lower coordinates therefore gives
\[
 \sum_{r,s,h\geq0}H_{(1^r,h+2,1^s)}X^rY^sZ^h
 =\frac1{(1-cX)(1-cY)(1-Z)}
 +\frac{t\mathcal R_{\Box}}{(1-X)(1-Y)(1-Z)}.
\]
Substitution of the preceding expression for $\mathcal R_{\Box}$ proves
\eqref{eq:singleton-sandwich-series}.

For the gamma transform, substitute $X=U/c$, $Y=V/c$, $Z=W/c$, and
$z=t/c^2$ in \eqref{eq:singleton-sandwich-series}, and divide the result by
$c$.  Simplification gives \eqref{eq:singleton-sandwich-gamma-series}.  Its
coefficientwise nonnegativity is transparent after putting
\[
 a=\frac{U}{1-U},\qquad
 d=\frac{V}{1-V},\qquad
 e=\frac{W}{1-W}.
\]
Indeed, the denominator in \eqref{eq:singleton-sandwich-gamma-series} is
\[
 (1-U)(1-V)(1-W)\bigl(1-z(ad+ae+de+ade)\bigr),
\]
and hence the right-hand side equals
\[
 \frac{1+z(a+d)}{(1-U)(1-V)(1-W)}
 \sum_{k\geq0}z^k(ad+ae+de+ade)^k.
\]
Every coefficient is nonnegative, which proves gamma-positivity.
\end{proof}

Setting $V=0$ recovers the $b\geq2$ part of the rectangular series in
Theorem~\ref{thm:singleton-prefix-rectangles}.  The new feature is that the
universal denominator has advanced from $D_2$ to $D_3$, even though the
number of layers in the host is still unbounded.  Determining whether this
larger gamma-positive family is also real-rooted and interlacing is the next
natural problem.


\subsection{A marked-slot model at fixed central size}
\label{subsec:sandwich-slots}

Fix the central layer size $b\geq2$, and put $m=b-1$. The $m$ slots below
may be indexed by the internal bonds of that layer; the two ordered sets
of flank positions are $[r]$ and $[s]$.

\begin{proposition}[Fixed-central-layer coefficients]
\label{prop:sandwich-slot-model}
Use the convention $\binom nk=0$ for $k<0$ or $k>n$.
For $m\geq1$, the gamma polynomial of $H_{(1^r,m+1,1^s)}$ is
\begin{equation}
 \Delta_{r,s,m-1}(z)
 =\sum_{\substack{i,j,k\geq0\\i+j\leq m}}
 \binom{m}{i,j,m-i-j}\binom{m+k-1}{k}
 \binom{r}{k+i}\binom{s}{k+j}z^{k+i+j}.
 \label{eq:sandwich-slot-coefficients}
\end{equation}
In particular, its degree is exactly
\begin{equation}
 \min\left\{r+s,\ r+m,\ s+m,\
              \left\lfloor\frac{r+s+m}{2}\right\rfloor\right\}.
 \label{eq:sandwich-exact-gamma-degree}
\end{equation}
\end{proposition}

\begin{proof}
Let $a=U/(1-U)$ and $c=V/(1-V)$. Extracting $W^{m-1}$ from
\eqref{eq:singleton-sandwich-gamma-series} gives
\begin{equation}
 \sum_{r,s\geq0}\Delta_{r,s,m-1}(z)U^rV^s
 =\frac1{(1-U)(1-V)}
   \left(\frac{1+z(a+c)}{1-zac}\right)^m.
 \label{eq:sandwich-central-power}
\end{equation}
Indeed, the denominator of that series, viewed as a polynomial in $W$,
is $(1-U)(1-V)\bigl((1-zac)-W(1+z(a+c))\bigr)$.
The multinomial and negative-binomial expansions, together with
$[U^r]a^u/(1-U)=\binom ru$, prove
\eqref{eq:sandwich-slot-coefficients}.

Here is a direct interpretation of its summands. Choose $k+i$ left and
$k+j$ right flank positions. Choose disjoint sets of $i$ left-marked and
$j$ right-marked slots among the $m$ slots, and a weak composition
$(k_1,\ldots,k_m)$ of $k$. Read the slots in order. At a marked slot,
first take the next selected position on its specified side; then take
$k_\ell$ pairs, each consisting of the next selected position on each side.
Thus every selected position is used once, and the record has $k+i+j$
marks, counting each pair as one mark. These records are enumerated by
\eqref{eq:sandwich-slot-coefficients}. This describes the gamma coefficients
through ordered flank positions and central slots; it does not identify
the records with individual chains in the interval.

For the degree, put $v=k+i+j$. Feasibility gives
$v\leq r+s$, $v\leq r+m$, $v\leq s+m$, and $2v\leq r+s+m$.
Conversely, if $v$ is the minimum in
\eqref{eq:sandwich-exact-gamma-degree}, take
\[
 k=\max(0,v-m),\qquad j=\max(0,v-r),\qquad i=v-k-j.
\]
The four bounds imply $i\geq0$, $k+i\leq r$, $k+j\leq s$,
and $i+j\leq m$. The corresponding summand is positive, proving
attainability.
\end{proof}

For instance, $[z]\Delta_{r,s,m-1}=m(r+s+rs)$: a record with one mark
chooses one of the $m$ slots and either one left position, one right
position, or one pair of positions.

The same formula provides a finite sum for every central size, without
expanding the full three-variable series. Its first two values of $m$
will also allow exact comparisons of zeros.


\subsection{The central layer of size two}
\label{subsec:central-two}

The first nonrectangular slice of the singleton-sandwich family is still
completely real-rooted.  Its proof uses the rectangular Jacobi polynomial as
an interlacer.

\begin{theorem}[Central-two real-rootedness]
\label{thm:central-two-real-rootedness}
For $r,s\geq0$, let
\[
 J_{r,s}(z)=\sum_{k\geq0}\binom{r}{k}\binom{s}{k}z^k.
\]
The gamma polynomial $\Delta_{r,s,0}(z)$ of
$H_{(1^r,2,1^s)}(t)$ satisfies
\begin{equation}
\label{eq:central-two-integral}
 \Delta_{r,s,0}(z)
 =(1-2z)J_{r,s}(z)
 +(r+s+2)\int_0^z J_{r,s}(u)\,du.
\end{equation}
Unless $r=s=0$, it has only simple negative zeros, and its zeros strictly
interlace those of $J_{r,s}$, with the constant--linear convention when one
parameter is zero.  Consequently, $H_{(1^r,2,1^s)}(t)$ is real-rooted and
\[
 H_{(1^r,s+1)}(t)\preceq H_{(1^r,2,1^s)}(t).
\]
\end{theorem}

\begin{proof}
Set $W=0$ in \eqref{eq:singleton-sandwich-gamma-series}.  Since
\[
 \frac1{1-U-V+(1-z)UV}
 =\sum_{r,s\geq0}J_{r,s}(z)U^rV^s,
\]
coefficient extraction gives
\begin{align*}
 \Delta_{r,s,0}(z)
 ={}&J_{r,s}(z)
 +z\sum_{i=0}^{r-1}J_{i,s}(z)
 +z\sum_{j=0}^{s-1}J_{r,j}(z)\\
 ={}&J_{r,s}(z)
 +\sum_{k\geq0}
 \left(
  \binom{r}{k+1}\binom{s}{k}
  +\binom{s}{k+1}\binom{r}{k}
 \right)z^{k+1}.
\end{align*}
Writing $I(z)=\int_0^zJ_{r,s}(u)\,du$, the two $z$-weighted sums on the first line are
respectively $(r+1)I-zJ_{r,s}$ and $(s+1)I-zJ_{r,s}$.  This proves
\eqref{eq:central-two-integral}.

By symmetry, assume $r\leq s$.  The case $r=0$ is immediate from
$\Delta_{0,s,0}(z)=1+sz$, so suppose $r\geq1$.  The Jacobi representation
\eqref{eq:delannoy-jacobi} shows that $J_{r,s}$ has $r$ simple negative
zeros.  A coefficient comparison also gives
\begin{equation}
\label{eq:central-two-integral-elimination}
 (r+1)(s+1)I(z)
 =z(1-z)J'_{r,s}(z)+(r+s+1)zJ_{r,s}(z).
\end{equation}
At any zero $\rho$ of $J_{r,s}$, equations
\eqref{eq:central-two-integral} and
\eqref{eq:central-two-integral-elimination} therefore give
\[
 \Delta_{r,s,0}(\rho)
 =\frac{(r+s+2)\rho(1-\rho)}{(r+1)(s+1)}
  J'_{r,s}(\rho).
\]
The prefactor is negative.  Hence the signs of
$\Delta_{r,s,0}$ alternate at consecutive zeros of $J_{r,s}$, giving one
zero in every intervening interval.  At the largest zero, the displayed
value is negative, whereas $\Delta_{r,s,0}(0)=1$, so there is one further
zero before $0$.

If $s=r$, these are already $r=\deg\Delta_{r,r,0}$ distinct negative
zeros.  If $s>r$, then $\deg\Delta_{r,s,0}=r+1$ with positive leading
coefficient.  At the smallest zero of $J_{r,s}$ the sign of
$\Delta_{r,s,0}$ is $(-1)^r$, while its sign at negative infinity is
$(-1)^{r+1}$; this supplies the final zero to the left.  Thus all zeros are
simple and the asserted strict interlacing holds.  The gamma--interlacing
correspondence now gives the weak Chow interlacing, because $J_{r,s}$ is the
gamma polynomial of $H_{(1^r,s+1)}(t)$ by
Theorem~\ref{thm:singleton-prefix-rectangles}.
\end{proof}


\subsection{The central layer of size three}
\label{subsec:central-three}

The first central growth step also admits a uniform theorem, with no
restriction on the two flank lengths.

\begin{theorem}[Central-three interlacing]
\label{thm:central-three-interlacing}
For all $r,s\geq0$, the polynomial $H_{(1^r,3,1^s)}(t)$ has only negative
real zeros, and
\begin{equation}
 H_{(1^r,2,1^s)}(t)\preceq H_{(1^r,3,1^s)}(t).
 \label{eq:central-two-three-interlacing}
\end{equation}
Unless $r=s=0$, the gamma polynomials $\Delta_{r,s,0}$ and
$\Delta_{r,s,1}$ have simple negative zeros and strictly interlace.
\end{theorem}

\begin{proof}
By symmetry assume $r\leq s$. The case $r=0$ follows from
Theorem~\ref{thm:singleton-prefix-rectangles}; at $r=s=0$ the Chow
polynomials are $1+t$ and $(1+t)^2$. For the remaining cases write
\[
 P=\Delta_{r,s,0},\qquad Q=\Delta_{r,s,1},\qquad J=J_{r,s}.
\]
Their coefficients are given by \eqref{eq:central-two-integral} and
\eqref{eq:sandwich-slot-coefficients}. In particular,
\begin{align}
 [z^j]Q={}&(j+1)\binom rj\binom sj
 +2j\left(\binom rj\binom{s}{j-1}
                 +\binom{r}{j-1}\binom sj\right)\notag\\
 &+(j-1)\left(\binom rj\binom{s}{j-2}
       +2\binom{r}{j-1}\binom{s}{j-1}
       +\binom{r}{j-2}\binom sj\right).
 \label{eq:central-three-coefficients}
\end{align}

We first treat $r=1<s$. In this case
\[
 P=1+(2s+1)z+\binom s2z^2,\qquad
 Q=1+(4s+2)z+\left(5\binom s2+2s\right)z^2+2\binom s3z^3.
\]
Let $w_s=4s^3+25s^2+35s+8$. The remainder of $Q$ on division by $P$ is
\[
 -\frac{w_sz+4(s+1)(s+2)}{3s(s-1)}.
\]
Its zero $\eta=-4(s+1)(s+2)/w_s$ satisfies
\[
 P(\eta)=-\frac{s(s-1)(8s^4+80s^3+211s^2+181s+24)}{w_s^2}<0.
\]
Thus $\eta$ lies between the two negative zeros of $P$, and $Q$ has
respectively positive and negative signs there. Since $Q(0)=1$ and its
leading coefficient is positive, these signs give the asserted negative
zeros and interlacing: $Q$ has degree two for $s=2$ and degree three
for $s\geq3$.

We now treat $r=s\geq1$ and $2\leq r<s$. Put
\[
 S=r+s,\quad p=rs,\quad d=s-r,\quad K=(r+1)(s+1),
 \qquad A=4p+3S+2,\quad C=2p+3S+5,
\]
and define
\begin{align*}
 \mathcal R(x)&=K+Ax-d^2x^2,\\
 E&=Kd^2+(4-p)S+6,\qquad
 L=4p^2+8pS+5p+3S^2-5,\\
 \mathcal N(x)&=Cx\mathcal R(x)-Ex^2+Lx+K(K-1).
\end{align*}
At every zero $\rho$ of $P$, coefficient comparison gives
\begin{equation}
 Q(\rho)=-\frac{\mathcal N(-\rho)}
 {(S+2)(r+2)(s+2)(1-\rho)}J(\rho).
 \label{eq:central-three-sign-identity}
\end{equation}
For explicit verification, the polynomial identity behind this comparison is
\[
 (S+2)(r+2)(s+2)(1-z)Q(z)
 =T(z)P(z)-\mathcal N(-z)J(z),
\]
where
\begin{align*}
 T(z)={}&p^2+3pS+3S^2+3p+9S+8\\
 &+\bigl((p+S)d^2-pS-S^2-6p-8S-8\bigr)z
   +K(S+4-d^2)z^2.
\end{align*}
For a coefficientwise verification that also covers degree drops, put
$a_j=\binom rj\binom sj$ for $j\geq0$ and $a_j=0$ for $j<0$.
The constant coefficients of $J,P,Q$ are all $1$, and, for $j\geq1$,
\begin{align*}
 [z^j]P={}&a_j+\frac{S+2-2j}{j}a_{j-1},\\
 [z^j]Q={}&(j+1)a_j+2(S-j+1)a_{j-1}\\
 &+\frac{(r-j+1)(r-j+2)+(s-j+1)(s-j+2)}{j}a_{j-2}.
\end{align*}
These follow from \eqref{eq:central-two-integral} and
\eqref{eq:central-three-coefficients}. Substitution, using
$j^2a_j=(r-j+1)(s-j+1)a_{j-1}$, gives the stated polynomial
identity. Only positive integer $j$ occurs in a denominator, so no
division by a vanishing binomial coefficient is required at the
boundary; the constant coefficient is checked separately.

We need $\mathcal N(-\rho)>0$. If $r=s=n\geq1$, direct factorization gives
\[
 \mathcal N(x)=(n+1)(n+2)
 \bigl(2(4n^2+7n+1)x^2+2n(3n+4)x+n(n+1)\bigr)>0
 \quad(x>0).
\]
Suppose henceforth that $2\leq r<s$. The root orders in the proof of
Theorem~\ref{thm:central-two-real-rootedness}, together with the
positive leading coefficients, give $P'(\rho)/J(\rho)>0$; in
particular, $J(\rho)\ne0$.
To turn this into a useful bound, write, by
\eqref{eq:central-two-integral}--\eqref{eq:central-two-integral-elimination},
\[
 \begin{gathered}
 P(z)=A_0(z)J(z)+B_0(z)J'(z),\\
 A_0(z)=1+\left(\frac{(S+2)(S+1)}K-2\right)z,\qquad
 B_0(z)=\frac{S+2}{K}z(1-z).
 \end{gathered}
\]
The hypergeometric differential equation
\[
 z(1-z)J''(z)+\bigl(1+(S-1)z\bigr)J'(z)-pJ(z)=0
\]
eliminates $J''$ after differentiating $P=A_0J+B_0J'$.
At $P(\rho)=0$, substitute $J'(\rho)=-A_0(\rho)J(\rho)/B_0(\rho)$
to obtain
\begin{equation}
 -B_0(\rho)\frac{P'(\rho)}{J(\rho)}
 =\frac{\mathcal R(-\rho)}K.
 \label{eq:central-two-root-constraint}
\end{equation}
Since $B_0(\rho)<0$, we have $\mathcal R(-\rho)>0$.

In fact, $\mathcal N(x)>0$ whenever $x>0$ and $\mathcal R(x)>0$
in this parameter range. If $E\leq0$, this is immediate from its
definition, because $C,L,K(K-1)>0$. If $E>0$, put
\[
 U=(p-4)S-6,\qquad V=(p-5)S-p-7.
\]
Polynomial division gives
\[
 \mathcal N(x)=\left(Cx+\frac{E}{d^2}\right)\mathcal R(x)
 +\frac{Vd^2+UA}{d^2}x+\frac{K(U-d^2)}{d^2}.
\]
Here $p\geq2S-4$ and $S\geq5$, so
\[
 U-d^2=(S+4)p-S^2-4S-6\geq S^2-22>0.
\]
For $S\geq6$ we also have
$V\geq2S^2-11S-3>0$, making every displayed term positive.
The sole remaining pair $(r,s)=(2,3)$ has $Vd^2+UA=156>0$.
This proves \eqref{eq:central-three-sign-identity} has a strictly negative
prefactor at every zero of $P$.

It remains to count the zeros and record their order.
For $d=0$, the degrees of $(J,P,Q)$ are $(r,r,r+1)$.
For $d=1$ they are $(r,r+1,r+1)$, and for $d\geq2$ they are
$(r,r+1,r+2)$, by \eqref{eq:sandwich-exact-gamma-degree}.
When $d=0$, the known order is
$\alpha_1<\rho_1<\cdots<\alpha_r<\rho_r<0$ for the zeros of $(J,P)$.
When $d>0$, there is also a zero of $P$ to the left of $\alpha_1$.
At the largest zero of $P$, $J$ is positive, so $Q$ is negative,
whereas $Q(0)=1$.
Equation~\eqref{eq:central-three-sign-identity} therefore gives one zero
of $Q$ between every two consecutive zeros of $P$, and one between the
largest zero of $P$ and $0$. When $\deg Q=\deg P+1$, the signs at negative
infinity supply one additional zero to the left. These account for the
full degree, so all zeros are simple and negative.
In the equal-degree case the zeros of $Q$ lie to the right of the
corresponding zeros of $P$; in the other case they surround those of $P$.
Lemma~\ref{lem:gamma-chow-interlacing} now gives
\eqref{eq:central-two-three-interlacing}, including the weak alternation
at the common root $-1$.
\end{proof}


\subsection{Long central layers and a Laguerre differential transform}
\label{subsec:sandwich-asymptotics}

The full sandwich family also admits an all-size asymptotic statement:
for each fixed pair of flanks, sufficiently long central layers give
real-rooted Chow polynomials and consecutive interlacing.

\begin{theorem}[Eventual sandwich interlacing]
\label{thm:eventual-sandwich}
Fix $r,s\geq0$, and put $R=r+s$. For $0\leq j\leq R$, the coefficient
$[z^j]\Delta_{r,s,h}(z)$ is a polynomial in $h\geq0$ of degree at most $j$,
and $\deg\Delta_{r,s,h}\leq R$. Moreover,
\begin{equation}
 \Delta_{r,s,h}(x/h)\longrightarrow
 P_{r,s}(x):=\frac1{s!}\left(1+\frac{d}{dx}\right)^s
                 \bigl(x^s L_r(-x)\bigr)
 \label{eq:sandwich-laguerre-limit}
\end{equation}
coefficientwise as $h\to\infty$, where $L_r$ is the ordinary Laguerre
polynomial. The polynomial $P_{r,s}=P_{s,r}$ has degree $R$, constant
term $1$, and, if $R>0$, only simple negative zeros.
There is an integer $B_{r,s}\geq3$ such that, for all $b\geq B_{r,s}$,
\[
 H_{(1^r,b-1,1^s)}(t)\preceq H_{(1^r,b,1^s)}(t),
\]
and both polynomials have only negative real zeros.
\end{theorem}

\begin{proof}
Set $d=(1-U)(1-V)$, $D=d-zUV$, $T=U+V-UV$, and
$N=1+zU/(1-U)+zV/(1-V)$. Extracting the coefficient of $W^h$ in
\eqref{eq:singleton-sandwich-gamma-series} gives
\begin{equation}
 \sum_{r,s\geq0}\Delta_{r,s,h}(z)U^rV^s
 =\frac ND\left(1+\frac{zT}{D}\right)^h.
 \label{eq:sandwich-fixed-flanks-transfer}
\end{equation}
For a fixed coefficient $U^rV^s$, the binomial expansion on the right
terminates at $k=R$, since $T^k$ has total degree at least $k$.
Its $k$th summand has a factor $\binom hk z^k$, and its remaining
factors expand in nonnegative powers of $z$. Thus only $k\leq j$
can contribute to $[z^j]$, proving polynomial dependence on $h$ of
degree at most $j$. The nonnegative expansion used in the
proof of Theorem~\ref{thm:singleton-sandwich-series} also shows that each
factor contributing a power of $z$ uses at least one $U$ or $V$;
thus the degree in $z$ is at most $R$.

Substitute $z=x/h$ in \eqref{eq:sandwich-fixed-flanks-transfer}.
In the quotient by $(U^{r+1},V^{s+1})$ all relevant expansions are finite,
so the limit can be taken coefficientwise. It is
\[
 P_{r,s}(x)=[U^rV^s]\frac1d\exp\left(x\frac{T}{d}\right).
\]
Writing $a=U/(1-U)$ and $c=V/(1-V)$ gives $T/d=a+c+ac$ and makes
symmetry in $r,s$ explicit. The Laguerre generating function yields
\[
 P_{r,s}(x)=\sum_{i=0}^r\binom ri\frac{x^i}{i!}L_s^{(i)}(-x).
\]
Here $L_s^{(i)}$ denotes the generalized Laguerre polynomial with
parameter $i$. Rodrigues' formula,
$L_s^{(i)}(-x)=x^{-i}(1+d/dx)^s x^{s+i}/s!$,
then gives \eqref{eq:sandwich-laguerre-limit}; see
\cite[Chapter V]{Szego1939} for these classical identities.

To locate the zeros, start with $x^sL_r(-x)$. Its nonzero zeros are simple
and negative, and its zero at $0$ has multiplicity $s$.
For any nonconstant real-rooted polynomial $f$ with all zeros nonpositive,
$f+f'$ retains each multiple zero with multiplicity reduced by one.
Its other zeros are simple and lie strictly to the left of the largest
zero of $f$: this follows by solving
$f'/f=\sum_\rho m_\rho/(x-\rho)=-1$, whose left-hand side is
strictly decreasing between consecutive distinct zeros.
There is one solution to the left of the smallest zero and one between
each consecutive pair, and none to the right of the largest zero.
Applying this fact $s$ times gives $R$ simple negative zeros, unless
$R=0$, when $P_{0,0}=1$. The operator preserves degree, and its value
at $0$ after $s$ applications is $s!$, proving the remaining assertions.

Extend $\Delta_{r,s,h}$ to real $h$ through its coefficient polynomials.
If $R>0$, the polynomial coefficient bounds and the simple negative zeros of
$P_{r,s}$ now verify the hypotheses of
Lemma~\ref{lem:analytic-scaling-interlacing}. For all sufficiently
large real $h$, it gives all $R$ negative roots and the directed order
$\rho_i(h-1)<\rho_i(h)<\rho_{i+1}(h-1)$ for $i<R$, with
$\rho_R(h-1)<\rho_R(h)<0$. In particular, the derivative control
comes from analytic branches, not from differentiating a bare
asymptotic error term.
Since $h=b-2$, choose $B_{r,s}\geq\max(3,R+2)$ so that this order
holds and $2R\leq R+b-2$. Lemma~\ref{lem:gamma-chow-interlacing},
with $n=R+b-2$, transfers it to the asserted Chow interlacing.
For $R=0$ the claim follows directly from $H_{(b)}=(1+t)^{b-1}$.
\end{proof}


\section{Discussion and open problems}
\label{sec:discussion}

The fixed-layer theorem and the two unbounded-layer families indicate that
principal composition intervals are more rigid than arbitrary graded
posets, even though their order relations are not determined by rank alone.
The rational denominator \(D_k\) captures the common finite-dimensional
transfer, while the Delannoy--Jacobi formula shows that special
unbounded-dimensional slices can carry classical orthogonal-polynomial
structure.  Neither feature is visible from the general definition of the
characteristic Chow function.

The rank-uniformity required by the TN-poset theorem and its dual version
in \cite[Theorem~6.1]{BrandenVecchi2025TN} can fail in both orientations.
For example, in $P_{(1,2,1)}$, the rank-three elements $(1,1,1)$ and $(1,2)$
have respectively one and two rank-two elements below them. The rank-two
elements $(1,1)$ and $(2)$ have respectively three and two rank-three
elements above them. Thus that theorem does not apply directly to the
whole family; a different matrix realization is not ruled out.

The first unresolved question is the uniform two-layer statement.  The
proved range \(\min(a,b)\leq4\), the exact fixed-width transfer, and the
Laguerre asymptotics all support the expectation that
\[
 H_{(a,b-1)}(t)\preceq H_{(a,b)}(t)
 \qquad(a,b\geq1).
\]
Conjecture~\ref{conj:gamma-difference-sturm} is designed to be stronger than
this assertion and to reduce it to a concrete boundary family.  A proof by
an interlacing-preserving operator or by a totally nonnegative production
matrix would also connect the two-layer kernel directly with recent matrix
methods for Chow polynomials \cite{BrandenVecchi2025TN}.

The singleton-sandwich family is gamma-positive, with an explicit
marked-slot interpretation of its coefficients.
Theorems~\ref{thm:central-two-real-rootedness}
and~\ref{thm:central-three-interlacing} settle central sizes two and three
for all flank lengths, including the first central growth step.
Theorem~\ref{thm:eventual-sandwich} also settles sufficiently long central
layers for every fixed pair $(r,s)$. The remaining question is whether
every \(H_{(1^r,b,1^s)}(t)\), \(b\geq4\), is real-rooted, and whether each
growth step from \(b-1\) to \(b\) preserves interlacing.
The fixed-central-layer power formula
\eqref{eq:sandwich-central-power} offers a concrete starting point.
A multivariate stable refinement or an interlacing-preserving recurrence
would go beyond coefficientwise positivity and would parallel recent
stable refinements of Narayana polynomials \cite{YangZhang2025}.

Ultimately one may ask which compositions \(\beta\) have real-rooted
\(H_\beta(t)\), or whether real-rootedness holds for every principal
composition interval.  The topology of general layered-permutation
intervals is known to be delicate
\cite{McNamaraSteingrimsson2015}, so a positive answer would not be a formal
consequence of a uniform shellability assertion.  This makes principal
composition intervals a useful family in which to compare order-theoretic,
matrix, and analytic mechanisms for Chow real-rootedness.

The accompanying supplementary verification programs check coefficient and
symbolic identities and perform finite exact and numerical root audits.
They serve as consistency checks. The proofs establish the theorems for
all stated parameters and use no finite root audit as a premise.


\section*{Statements and declarations}

\paragraph{Funding.}
The author received no specific funding for this work.

\paragraph{Competing interests.}
The author declares no competing interests.

\paragraph{Data and supplementary computational materials.}
This theoretical study uses no external datasets. The formulas and proofs
supporting the results are contained in the article. The ancillary material
supplied with the preprint source provides supplementary verification scripts, dependency and execution
instructions, and reports specifying the scope of the checks. These
computations are consistency checks; finite tests are not premises of the
theorems.

\paragraph{Supplementary information.}
Supplementary verification programs, scope reports, a README, and
execution instructions are supplied with the preprint source in the
directory \path{anc/verification_code/}.

\paragraph{AI Disclosure.}
ChatGPT/Codex assisted with algebraic simplification, checking and refining
proof arguments, preparation of supplementary verification code, and
manuscript revision throughout the paper. This assistance materially
influenced mathematical derivations and their presentation, as well as the
accompanying code. The author retains full responsibility for the
correctness, originality, attribution, and presentation of all content.

\end{document}